\documentclass[aos,preprint]{imsart}
\usepackage[french,english]{babel}
\usepackage[latin1]{inputenc}
\usepackage{amsmath}
\usepackage{amssymb}
\usepackage{graphicx,epsfig}
\usepackage{bbm}
\usepackage[latin1]{inputenc}
\input epsf

 \numberwithin{equation}{section}
\allowdisplaybreaks 

\def\BibTeX{{\rm B\kern-.05em{\sc i\kern-.025em b}\kern-.08em
    T\kern-.1667em\lower.7ex\hbox{E}\kern-.125emX}}

\def\limitepsn{\renewcommand{\arraystretch}{0.5}
\begin{array}[t]{c}\stackrel{a.s.}{\longrightarrow} \\
{\scriptstyle n \rightarrow
\infty}\end{array}\renewcommand{\arraystretch}{1}}

\def\limitepsu{\renewcommand{\arraystretch}{0.5}
\begin{array}[t]{c}\stackrel{a.s.}{\longrightarrow} \\
{\scriptstyle ~ n,\delta \rightarrow
\infty}\end{array}\renewcommand{\arraystretch}{1}}

\def\limiteproban{\renewcommand{\arraystretch}{0.5}
\begin{array}[t]{c}\stackrel{{\cal P}}{\longrightarrow} \\
{\scriptstyle n\rightarrow
\infty}\end{array}\renewcommand{\arraystretch}{1}}

\newcommand{\be}{\begin{equation}}
\newcommand{\ee}{\end{equation}}
\newcommand{\bd}{\begin{displaymath}}
\newcommand{\ed}{\end{displaymath}}

\newcommand{\ba}{\begin{eqnarray}}
\newcommand{\ea}{\end{eqnarray}}
\newcommand{\ban}{\begin{eqnarray*}}
\newcommand{\ean}{\end{eqnarray*}}

\newcommand{\R} {I\!\!R}
\newcommand{\E} {I\!\! E}
\newcommand{\N} {I\!\! N}
\newcommand{\Z} {\mathbb{Z}}

\renewcommand{\arraystretch}{.8}

\renewcommand{\Box}{\hfill\rule{0.25cm}{0.25cm}} 

\newtheorem{Prop}{Proposition}[section]
\newtheorem{lem}{Lemma}[section]
\newtheorem{Theo}{Theorem}[section]
\newtheorem{cor}{Corollary}[section]
\newtheorem{rem}{Remark}[section]
\newtheorem{Ex}{Example}
\newenvironment{dem}{\  {\bf Proof }}
{\Box\par\medskip\noindent}

\def\1{{\bf 1}}

\begin{document}

\begin{frontmatter}

\title{Detecting multiple change-points in general causal time series using penalized quasi-likelihood}
\runtitle{Detecting multiple change-points in general causal time
series}
\begin{aug}
  \author{Jean-Marc Bardet \ead[label=e1]{Jean-Marc.Bardet@univ-paris1.fr}
  \ead[label=u1,url]{http://samm.univ-paris1.fr/-Jean-Marc-Bardet}}
  \author{William Kengne \ead[label=e2]{William.Kengne@malix.univ-paris1.fr}
  \ead[label=u3,url]{http://samm.univ-paris1.fr/William-Kengne}}
  \and
  \author{Olivier Wintenberger \ead[label=e3]{Olivier.Wintenberger@ceremade.fr}
  \ead[label=u2,url]{http://wintenberger.fr}}
  \runauthor{J.-M. Bardet, W. Kengne and O. Wintenberger}

  \affiliation{SAMM, Université Paris 1 Panth\'eon-Sorbonne, France \\
  CEREMADE, Université Paris-Dauphine, FRANCE}

\address{
SAMM, Université Paris 1\\
90 Rue de Tolbiac\\
75634 Paris Cedex 13, FRANCE\\
\printead{e1}\\
          \printead{u1}\\
          ~\\
          \printead{e2}\\
          \printead{u3}\\
          ~\\
          CEREMADE, Université Paris-Dauphine\\
Place du Maréchal de Lattre de Tassigny\\
75016 Paris, FRANCE\\
          \printead{e3}\\
          \printead{u2}}

\end{aug}

\begin{abstract}~~This paper is devoted to the off-line multiple change-point detection in a semiparametric framework. The time series is supposed to belong to a large class of models including AR($\infty$), ARCH($\infty$), TARCH($\infty$),... models where the coefficients change at each instant of breaks. The different unknown parameters (number of changes, change dates and parameters of successive models) are estimated using a penalized contrast built on  conditional quasi-likelihood. Under Lipshitzian conditions on the model, the consistency of the estimator is proved when the moment order $r$ of the process satisfies $r\geq 2$. If $r\geq 4$, the same convergence rates for the estimators than in the case of independent random variables are obtained. The particular cases of AR($\infty$), ARCH($\infty$) and TARCH($\infty$) show that our method notably improves the existing results.\\
\end{abstract}



\end{frontmatter}

\section{Introduction}
The problem of the detection of change-points is a classical problem as well as in the statistic than in the signal processing community. If the first important result in this topic was obtained by Page \cite{Page1955} in $1955$, real advances have been done in the seventies, notably with the results of Hinkley (see for instance Hinkley \cite{Hinkley1970}
) and the topic of change detection became a distinct and important field of the statistic since the eighties (see the book of Basseville and Nikiforov \cite{Basseville1993} for a large overview). \\
Two approaches are generally considered for solving a problem of change detection: an 'on-line' approach leading to sequential estimation and an 'off-line' approach which arises when the series of observations is complete. Concerning this last approach, numerous results were obtained for independent random variables in a parametric frame (see for instance Bai and Perron \cite{Bai1998}). The case of the off-line detection of multiple change-points in a parametric or semiparametric frame for dependent variables or time series also provided an important literature. The present paper is a new contribution to this problem. \\
~\\
In this paper, we consider a general class ${\cal M}_T(M,f)$ of
causal (non-anticipative) time series. Let $M$ and $f$ be a
measurable functions such that for all $(x_i)_{i\in \N}\in \R^{\N}$,
$M\big((x_i)_{i\in \N}\big)$ is a $(m\times p)$ non-zero real matrix
and $f\big((x_i)_{i\in \N}\big) \in \R^m$. Let $T \subset \Z$ and
$(\xi_t)_{t\in \mathbb{Z}}$ be a sequence of centered independent
and identically distributed (iid) $\mathbb{R}^p$-random vectors
called the innovations and satisfying $\textrm{var}(\xi_0)=I_p$ (the
identity matrix of dimension $p$). Then, define
~\\
{\bf Class ${\cal M}_T(M,f)$:} {\it The process $X=(X_t)_{t\in\Z}$
belongs to  ${\cal M}_T(M,f)$ if it satisfies the relation:}
\begin{equation}\label{model}
X_{t+1}=M \big ((X_{t-i})_{i\in \N}\big ) \xi_t+f\big
((X_{t-i})_{i\in \N}\big )\quad\mbox{for all $t\in T$}.
\end{equation}
\noindent The existence and properties of these general affine processes  were studied in Bardet and Wintenberger \cite{Bardet2009} as a particular case of chains with infinite memory considered in Doukhan and Wintenberger \cite{Doukhand}. Numerous classical real valued time series are included in ${\cal M}_\Z(M,f)$: for instance AR$(\infty)$, ARCH$(\infty)$, TARCH$(\infty)$, ARMA-GARCH or bilinear processes.  \\
~\\
The problem of change-point detection is the following: assume that
a trajectory $(X_1,\cdots,X_n)$ of $X=(X_t)_{t\in\Z}$ is observed
where
\begin{equation}\label{problem}
X \in {\cal M}_{T_j^*}(M_{\theta_j^*},f_{\theta_j^*})\quad \mbox{for
all $j=1,\ldots,K^*$},\qquad\mbox{with}
\end{equation}
\begin{itemize}
\item  $K^*\in \N^*$, $T_j^*=\{t^*_{j-1}+1,t^*_{j-1}+2,\ldots, t^*_{j}\}$ with $0<t_1^*<\ldots<t_{K^*-1}^*<n$, $t^*_j\in \N$ and by convention $t_0^*=-\infty$ and $t^*_{K^*}=\infty$;
\item $\theta^*_j=(\theta^*_{j,1},\cdots,\theta^*_{j,\, d})\in \Theta\subset \mathbb{R}^d$ for $j=1,\ldots,K^*$.
\end{itemize}
The aim in the problem is the estimation of the unknown  parameters $\big (K^*,(t^*_j)_{1\leq j \leq K^*-1}, (\theta^*_j)_{1\leq j \leq K^*}\big )$.  In the literature it is generally supposed that $X$ is a stationary process on each set $T^*_j$ and is independent on each $T^*_j$ from the other $T_k^*$, $k\neq j$ (for instance in \cite{LavielleMoulines2000},  \cite{Leipus2000}, \cite{David1995} and \cite{David2008}). Here the problem \eqref{problem} does not induce such assumption and thus the framework is closer to the applications, see Remark 1 in  \cite{David2008}. \\
~\\
In the problem of change-point detection,  numerous papers were devoted to the CUSUM procedure (see for instance Kokozska and Leipus \cite{Leipus2000}  in the specific case of ARCH$(\infty)$ processes).  In Lavielle and Ludena \cite{Lavielle2000} a "Whittle" contrast is used for estimating the break dates in the spectral density of piecewise long-memory processes (in a semi-parametric framework). Davis {\it et al.} \cite{David1995}  proposed a likelihood ratio as the estimator of break points for an AR$(p)$ process. Lavielle and Moulines  \cite{LavielleMoulines2000} consider a general contrast using the mean square errors for estimating the parameters. In Davis {\it et al.} \cite{David2008},  the criteria called Minimum Description Length (MDL) is applied to a large class of nonlinear time-series model.\\
We consider here a semiparametric estimator based on a penalized contrast using the quasi-likelihood function. For usual stationary time series,  the conditional quasi-likelihood
is constructed as follow:
\begin{enumerate}
\item Compute the conditional likelihood (with respect to $\sigma\{X_0,X_{-1},\ldots\}$) as if $(X_t)_{t\in \Z}$ is known and when the process of innovations is a Gaussian sequence;
\item Approximate this computation for a sample $(X_1,\ldots,X_n)$;
\item Apply this approximation even if the process of innovations is not a Gaussian sequence.
\end{enumerate}
The quasi-maximum likelihood estimator (QMLE) obtained by maximizing the quasi-likelihood function has convincing asymptotic properties in the case of GARCH processes (see Jeantheau \cite{Jeantheau1998},  Berkes {\it et al.} \cite{Berkes2003}, Franck and Zakoian \cite{Francq2004}) or generalizations of GARCH processes (see Mikosch and Straumann \cite{Straumann2006}, Robinson and Zaffaroni \cite{Robinson2006}). Bardet and Wintenberger \cite{Bardet2009} study the asymptotic normality of the QMLE of $\theta$ applied to  ${\cal M}_{\Z}(f_\theta,M_\theta)$. Thus, when $K^*$ is known, a natural estimator of the parameter $\big ((t^*_j)_{1\leq j \leq K^*-1}, (\theta^*_j)_{1\leq j \leq K^*}\big )$ for a process satisfying \eqref{problem} is the QMLE on every intervals $[t_j+1,\ldots,t_{j+1}]$ and every parameters $\theta_j$ for $1\le j\le K^*$. However we consider here that $K^*$ is unknown and such method cannot be directly used. The solution chosen is to penalize the contrast by an additional term $\beta_n K$, where $(\beta_n)_{n\in \N}$ is an increasing sequence of real numbers (see the final expression of the penalized contrast in \eqref{contrast}). Such procedure of penalization  was previously used for instance by Yao \cite{Yao1988} to estimate the number of change-points with the Schwarz criterion and by Lavielle and Moulines \cite{LavielleMoulines2000}. Hence the minimization of the penalized contrast leads to an estimator (see \eqref{estimateur}) of the parameters $\big (K^*, (t^*_j)_{1\leq j \leq K^*-1}, (\theta^*_j)_{1\leq j \leq K^*}\big )$. \\
Classical heuristics such as the BIC one or the MDL one of \cite{David2008} lead to choose $\beta_n \propto \log n$. In our study, such penalizations are excluded in some cases, when the models ${\cal M}_T(M,f)$ are very dependent of their whole past, see Remark \ref{ratevn} for more details. Finally, we will show that an ``optimal'' penalization is $\beta_n \propto \sqrt n$ which overpenalizes the number of breaks to avoid artificial breaks in cases of models very dependent of their whole past (see Remark \ref{choicevn}).\\
~\\
The main results of the paper are the following: under Lipshitzian condition on $f_\theta$ and $M_\theta$, the estimator $\big ( \widehat K_n,(\widehat t_j/n)_{1\leq j \leq \widehat K_n-1},(\widehat \theta_j)_{1\leq j \leq \widehat K_n} \big)$ is consistent when the moment of order $r$ on the innovations and $X$ is larger than $2$. If moreover Lipshitzian conditions are satisfied by the derivatives of $f_\theta$ and $M_\theta$ and if $r\geq 4$, then the convergence rate of $(\widehat t_j/n)_{1\leq j \leq \widehat K_n-1}$ is $O_P(n^{-1})$ and a Central Limit Theorem (CLT) for  $(\widehat \theta_j)_{1\leq j \leq \widehat K_n}$ (with a $\sqrt n$-convergence rate)  is established.  These results are "optimal" in the sense that they are the same than in an independent setting. \\
~\\
Section \ref{propmod} is devoted to the presentation
of the model and the assumptions and the study of the
existence of a nonstationary solution of the problem \eqref{problem}. The definition of the estimator and its
asymptotic properties are studied in Section \ref{estima}. The
particular examples of AR($\infty$), ARCH($\infty$) and
TARCH($\infty$) processes are detailed in Section \ref{Examples}. Section
\ref{proofs} contains the main proofs.

\section{Assumptions and existence of a solution of the change process}\label{propmod}

\subsection{Assumptions on the class of models ${\cal M}_\Z (f_\theta,M_\theta)$}
Let $\theta\in\R^d$ and $M_\theta$ and $f_\theta$ be numerical
functions such that for all $(x_i)_{i\in \N}\in \R^{\N}$,
$M_\theta\big((x_i)_{i\in \N}\big)\neq 0$ and
$f_\theta\big((x_i)_{i\in \N}\big) \in \R$. We use the following
different norms:
\begin{enumerate}
\item  $\|\cdot\|$ applied to a vector denotes the Euclidean norm of
    the vector;
    \item for any compact set $ \Theta\subseteq\R^d$ and for any
    $g:\Theta \longrightarrow\R^{d'}$;
    $
    \|g\|_\Theta=\sup_{\theta\in \Theta}(\|g(\theta)\|)$;
    \item  for all $x=(x_1,\cdots,x_K)\in \mathbb{R}^K,~~ \|x \|_m= \underset{i=1,\cdots,K} {\mbox{max}}|x_i|$;
    \item if $X$  is $\R^p$-random variable with  $r\geq 1$ order moment, we set $\|X \|_r = (\E\|X\|^r)^{1/r}$.
\end{enumerate}
Let  $\Psi_\theta=f_\theta, \, M_\theta$ and $i=0,\, 1, \,
2$, then for any compact set $ \Theta\subseteq\R^d$, define\\
~\\
{\bf Assumption A$_i(\Psi_\theta,\Theta)$}: {\em Assume that
$\|{\partial^i\Psi_\theta(0)}/{\partial\theta^i}\|_\Theta<\infty$
 and there exists a sequence of non-negative real number $(\alpha^{(k)}_i(\Psi_\theta,\Theta))_{i\geq 1}$ such that $ \sum\limits_{k=1}^{\infty} \alpha^{(i)}_k(\Psi_\theta,\Theta) <\infty$  satisfying
\begin{equation*}
\Big\|\dfrac{\partial^i\Psi_\theta(x)}{\partial\theta^i}-\dfrac{\partial^i\Psi_\theta(y)}{\partial\theta^i}\Big\|_\Theta
\leq \sum\limits_{k=1}^{\infty}\alpha^{(i)}_k(\Psi_\theta,\Theta)|x_k-y_k|\quad \mbox{for all}~x, y \in \R^{\N}.\\
\end{equation*}}
In the sequel we refer to the particular case called "ARCH-type process"
 if $f_\theta =0$ and if the following assumption holds on $
h_\theta = M_\theta^2$:\\
~\\
{\bf Assumption A$_i(h_\theta,\Theta)$}: {\em Assume that
$\|{\partial^i h_\theta(0)}/{\partial\theta^i}\|_\Theta<\infty$
 and there exists a sequence of non-negative real number $(\alpha^{(k)}_i(h_\theta,\Theta))_{i\geq 1}$ such as $ \sum\limits_{k=1}^{\infty} \alpha^{(i)}_k(h_\theta,\Theta) <\infty$  satisfying
\begin{equation*}
\Big\|\dfrac{\partial^i
h_\theta(x)}{\partial\theta^i}-\dfrac{\partial^i
h_\theta(y)}{\partial\theta^i}\Big\|_\Theta
\leq \sum\limits_{k=1}^{\infty}\alpha^{(i)}_k(h_\theta,\Theta)|x^2_k-y^2_k|\quad \mbox{for all}~x, y \in \R^{\N}.\\
\end{equation*}}
Now, for any $i=0, 1,2$ and $\theta\in \Theta$, under Assumptions
$A_i(f_\theta,\Theta)$ and $A_i(M_\theta,\Theta)$, denote:
\begin{eqnarray}\label{beta}
~~~~~~~~~\beta^{(i)}(\theta):=\sum_{k\geq
 1}\beta_{k}^{(i)}(\theta)\quad \mbox{where}\quad\beta_{k}^{(i)}(\theta):=\alpha^{(i)}_k(f_\theta,\{\theta\})+(\E |\xi_0|^r)^{1/r}
 \alpha^{(i)}_k(M_\theta,\{\theta\}),
\end{eqnarray}
and under Assumption $A_i(h_\theta,\Theta)$
\begin{eqnarray}\label{betatilde}
~~~~~~~~~\tilde\beta^{(i)}(\theta):=\sum_{k\geq
 1}\tilde\beta_{k}^{(i)}(\theta)\quad \mbox{where}\quad\tilde\beta_{k}^{(i)}(\theta):=(\E |\xi_0 |^r)^{2/r}
 \alpha^{(i)}_k(h_\theta,\{\theta\}).
\end{eqnarray}
The dependence with respect to $r$ of $\beta^{(i)}_k(\theta)$ and
$\tilde\beta^{(i)}_k(\theta)$ are omitted for notational
convenience. Then define:
\begin{eqnarray}\label{Thetar}
~~~\Theta(r):=\{\theta \in \Theta,\, A_0(f_\theta,\{\theta\})\, \,\mbox{and} \,A_0(M_\theta,\{\theta\}) \, \textrm{hold with}\, \beta^{(0)}(\theta)<1 \}&&\\
\nonumber&&\hspace{-8cm}\cup\{\theta \in \R^d,\,  f_\theta=0\,
\,\mbox{and} \,A_0(h_\theta,\{\theta\}) \, \textrm{holds with}\,
\tilde\beta^{(0)}(\theta)<1 \}.
\end{eqnarray}
From  \cite{Bardet2009} we have:
\begin{Prop}
If $\theta\in\Theta(r)$ for some $r\ge1$, there exists a unique
causal (non anticipative, i.e. $X_t$ is independent  of $(\xi_i)_{i>
t}$ for $t\in\Z$) solution $X=(X_t)_{t\in \Z}\in {\cal
M}_\Z(f_\theta,M_\theta)$ which is stationary, ergodic and satisfies
$ \|X_0\big\|_r <\infty$.
\end{Prop}
\begin{rem}
The Lipschitz-type hypothesis $A_i(\Psi_\theta,\Theta)$ are
classical when studying the existence of solution of general models.
For instance, Duflo \cite{Duflo1990}  used such a
Lipschitz-type inequality to show the existence  of Markov
chains. The subset $\Theta(r)$ is defined as a reunion to consider
accurately general causal models and ARCH-type models
simultaneously: for ARCH-type models $A_0(h_\theta,\{\theta\})$ is
less restrictive than $A_0(M_\theta,\{\theta\})$. However, remark
that $A_0(h_\theta,\{\theta\})$ is still not optimal for ensuring
the existence of a stationary solution for ARCH-type models.
\end{rem}
Let $\theta \in \Theta(r)$ and $X=(X_t)_{t\in \Z}$ a stationary solution included in ${\cal M}_\Z(f_\theta,M_\theta)$. For studying QMLE properties, it is convenient to assume the following assumptions:   \\
~\\
{\bf Assumption D$(\Theta)$:} $\exists\underline{h}>0$ such that
$\displaystyle \inf_{\theta \in
 \Theta}(|h_\theta(x)|)\geq \underline{h}$ for all $x\in \R^{\N}.$ \\
~\\
{\bf Assumption Id($\Theta$):} For all $\theta,\theta'\in \Theta^2$,
$$ \Big( f_{\theta}(X_0,X_{-1},\cdots)=f_{\theta'}(X_0,X_{-1},\cdots)~\mbox{and}~h_{\theta}(X_0,X_{-1},\cdots)=h_{\theta'}(X_0,X_{-1},\cdots) \ \text{a.s.}\Big) \Rightarrow \ \theta = \theta'.$$
{\bf Assumption Var($\Theta$):} For all $\theta  \in \Theta $,
$$ \big( \dfrac{\partial f_{\theta}}{\partial \theta^{(i)}}(X_0,X_{-1},\cdots) \big)_{1\leq i \leq d} \neq 0 \quad
\mbox{or} \quad \big( \dfrac{\partial h_{\theta}}{\partial
\theta^{(i)}}(X_0,X_{-1},\cdots) \big)_{1\leq i \leq d} \neq 0\quad
a.s.
$$
Assumption D$(\Theta)$ will be required to define the QMLE,
Id($\Theta$) to show the consistence of the QMLE and Var($\Theta$)
to show the asymptotic normality.

\subsection{Existence of the solution to the problem \eqref{problem}}
Consider the problem \eqref{problem} and let $(X_1,\ldots,X_n)$ be
an observed path of $X$. Then the past of $X$ before the time $t=0$
depends on $\theta_1^*$ and the future after $t=n$ depends on
$\theta_{K^*}^*$. The number $K^* -1$ of breaks, the instants $
t^*_1,\cdots,t^*_{k^*-1}$ of breaks and parameters
$\theta_1^*,\cdots,\theta_{K^*}^*$ are
unknown. Consider first the following notation. \\
~\\
\large{\textbf{Notation}}. {\em
\begin{itemize}
\item For $K \geq 2$,
    $
    \mathcal{F}_{K}=\big\{\underline{t}=(t_1,\ldots t_{K-1})\; ;\;0<t_1<\ldots<t_{K-1}<n\big\}.$
    In particular,
    $
    \underline{t}^*=\big(t_1^*,\ldots,t_{K^*-1}^*\big)\in \mathcal{F}_{K^*}
    $
    is the true vector of instants of change;

\item For $K \in \N^*$ and  $\underline{t}\in \mathcal{F}_{K} $, \,
     $
     T_k=\big\{t\in \mathbb{Z}
     $,
     $t_{k-1}<t\le t_{k}\big\}$ and $n_k=\textrm{Card}(T_k)$ with $1\leq k \leq K$. In particular;
     $
     T^*_j=\big\{t\in \mathbb{Z},\; t^*_{j-1}<t\le t^*_{j}\big\}
     $
     and $n_j^*=\textrm{Card}(T^*_j)$ for $1\leq j \leq K^*$. For all $1 \leq k \leq K$  and  $1 \leq j \leq K^*$,
     let
     $n_{kj}= \textrm{Card}(T^*_j \cap T_k)$;
 \end{itemize}}
\noindent The following proposition establishes the existence of the
nonstationary solution of the problem  \eqref{problem} and its
moments properties.
\begin{Prop} \label{prop1}
Consider the problem \eqref{problem}. Assume there exists $r\ge 1$
such that $\theta_j^*\in\Theta(r)$ for all $ j=1,\ldots,K^*$. Then
\begin{enumerate}
 \item[(i)] there exists a process $X=(X_t)_{t\in\Z}$ solution of the model \eqref{problem} such as $\|X_t\|_r<\infty$ for $t\in \Z$ and $X$ is a causal time series.
 \item[(ii)] there exists a constant $C>0$ such that for all  $t \in \Z$ we have $ \|X_t\|_r \le C$.
\end{enumerate}
\end{Prop}
\begin{rem}
The problem \eqref{problem} distinguishes the case $t\in T_1^*=\{1,\ldots,t^*_1\}$ to the other ones since it is easy to see that $(X_t)_{t\in T^*_1}$ is a stationary process while $(X_t)_{t> t^*_1}$ is not. However, all the results of this paper hold if $(X_t)_{t\in T^*_1}$ is defined as the other  $(X_t)_{t\in T^*_j}$, $j\geq 2$ (by defining a break in $t=0$) or, for instance, if we set $X_t=0$ for $t\leq 0$.
\end{rem}

\section{Asymptotic results of the estimation procedure}\label{estima}
\subsection{The estimation procedure}
The estimation procedure of the number of breaks $K^*-1$, the
instants of breaks $\underline t^*$ and the parameters $\underline
\theta^*$ is based on the minimum of a penalized contrast.
It is clear that if $(\xi_t)_t$ is a Gaussian process and if $X\in
{\cal M}_T(f_\theta,M_\theta)$ then for $s\in
T=\{t,t+1,\ldots,t'\}$, the distribution of $ X_{s}~|~(X_{s-j})_{j\in \N^*}$ is $
\mathcal{N}\big(f_{\theta}\big(X_{s-1},\ldots\big),h_{\theta}\big(X_{s-1},\ldots\big)\big)$.
Therefore, with the notation
 $ f^s_{\theta}=f_{\theta}\big(X_{s-1},X_{s-2}\ldots\big) $, $ M^s_{\theta}=M_{\theta}\big(X_{s-1},X_{s-2}\ldots\big) $ and
   $ h^s_{\theta}= {M^{s}_{\theta}}^2 $,
we deduce the conditional $\log$-likelihood on $T$ (up to an additional
constant) 
$$
L_n(T,\theta):=-\frac{1}{2}\sum\limits_{s\in T}
q_s(\theta)\;\;\textrm{with}\;\;q_s(\theta):=\frac{\left(X_s-f^s_{\theta}\right)^2}{h^s_{\theta}}+\log\left(h^s_{\theta}\right).
$$
By convention, we set $L_n(\emptyset, \theta_k):=0$. Since only
$X_1,\ldots,X_n$ are observed, $L_n(T,\theta)$ cannot be computed
because it depends on the past values $(X_{-j})_{j\in \N}$. We
approximate it by:
$$
\widehat L_n(T,\theta):=-\frac{1}{2}\sum\limits_{s\in T}\widehat
q_s(\theta)\quad
\textrm{where}\;\;\;\widehat{q}_s(\theta):=\frac{\big(X_s-\widehat{f}^s_{\theta}\big)^2}{\widehat{h}{^s_{\theta}}}
+\log\big(\widehat{h}{^s_{\theta}}\big)
$$
with
 $ \widehat{f}^t_{\theta}=f_{\theta}\big(X_{t-1},\ldots,X_{1},u\big)$, $ \widehat{M}^t_{\theta}=
 M_{\theta}\big(X_{t-1},\ldots,X_{1},u\big) $ and $\widehat{h}^t_{\theta}=( \widehat{M}^{t}_{\theta} )^2$ for any deterministic sequence $u=(u_n)$ with finitely many non-zero values.

\begin{rem}
For convenience, in the sequel we chose $u=(u_n)_{n\in
\mathbb{N}}$ with $u_n=0$ for all $n\in \mathbb{N}$ as in
\cite{Francq2004} or in \cite{Bardet2009}. Indeed, this choice has
no effect on the asymptotic behavior of estimators.
\end{rem}
Now, even if the process $(\xi_t)_t$ is non-Gaussian and for any
number of breaks $K-1\ge 1$ and any $\underline{t}\in
\mathcal{F}_{K},~\underline{\theta}\in \Theta(r)^{K}$, define the
contrast function $\widehat J_n$ by the expression:
\begin{eqnarray}\label{JN}
 \widehat J_n(K,\underline{t},\underline{\theta}):= -2\sum_{k=1}^{K}  \widehat L_n(T_k,\theta_k)= -2\sum_{k=1}^{K} \sum_{j=1}^{K^*} \widehat L_n(T_k\cap T^*_j,\theta_k),
\end{eqnarray}
~\\
Finally, let $(v_n)_{n\in \N}$ and $(\beta_n)_{n\in \N}$ be sequences satisfying $v_n\ge 1$ and $\beta_n:=n/{v_n} \to \infty$ ($n\to \infty$). Let   $K_{\max}\in \N^*$ and for $K \in \{1, \cdots, K_{\max} \}$ and $(\underline{t}, \underline{\theta}) \in \mathcal{F}_{K}
 \times\Theta(r)^K$ ($\Theta(r)$ is supposed to be a compact set) define the penalized contrast $\widetilde{J}_n$ by
\begin{eqnarray}\label{contrast}
\widetilde{J}_n ( K, \underline{t}, \underline{\theta} ): =
\widehat{J}_n (K, \underline{t}, \underline{\theta}) + \frac{n}{v_n}
\, K=\widehat{J}_n (K, \underline{t}, \underline{\theta}) + \beta_n
\, K
\end{eqnarray}
and the penalized contrast estimator $(\widehat{K}_n,
\widehat{\underline{t}}_n, \widehat{\underline{\theta}}_n)$ of
$(K^*,\underline{t}^*,\underline{\theta}^*)$ as
\begin{eqnarray}\label{estimateur}
(\widehat{K}_n, \widehat{\underline{t}}_n,
\widehat{\underline{\theta}}_n) = \underset{1\leq K\leq K_{\max}}
  {\mbox{Argmin}} ~~ \underset{( \underline{t}, \underline{\theta}) \in \mathcal{F}_{K} \times \Theta(r)^K } {\mbox{Argmin}}
  ( \widetilde{J}_n ( K, \underline{t}, \underline{\theta} ))  \quad \text{and} \quad \widehat{\underline{\tau}}_n =
  \dfrac{\widehat{\underline{t}}_n}{n}.
\end{eqnarray}
For $K\geq 1$ and $\underline{t} \in \mathcal{F}_{K}$, denote $
\widehat{\underline\theta}_n(\underline{t})=(\widehat \theta(T_1),\ldots,\widehat \theta(T_K)):=\underset{\underline{\theta}\in
\Theta(r)^{K}} {\mbox{Argmin}}\left(\widetilde
J_n(\underline{t},\underline{\theta})\right)=\underset{\underline{\theta}\in
\Theta(r)^{K}} {\mbox{Argmin}}\left(\widehat
J_n(\underline{t},\underline{\theta})\right)$. Then, $\widehat{\theta}_n(T_k)$ is the  QMLE of $\theta^{*}_{k}$ computed on $T_k$ and $\widehat{\underline\theta}_n(\underline{t}^*)$ is the $QMLE$ of $\underline \theta^*$.
\begin{rem}
If $K^\ast$ is known, the estimator of
$(\underline{\tau}^*,\underline{\theta}^*)$ may be obtained by
minimizing $\widehat {J}_n$ instead of $\widetilde{J}_n$. However
the knowledge of $K^\ast$ does not improve the asymptotic results
established in this paper.
\end{rem}

\subsection{Consistency of $(\widehat{K}_n, \widehat{\underline{t}}_n, \widehat{\underline{\theta}}_n)$}

\noindent For establishing the consistency, we add the couple of following classical assumptions in the problem of break detection:\\
~\\
{\bf Hypothesis B}: {\em  $\underset{j=1,\cdots,K^*-1} {\mbox{min}}\|\theta_{j+1}^*-\theta_j^*\|>0$.} \\
~\\
Furthermore, the distance between instants of breaks cannot be too small:\\
~\\
{\bf Hypothesis C}: {\em there exists $\tau_1^*,\ldots,\tau_{K^*-1}^*$ with $ 0<\tau_1^*<\ldots<\tau_{K^*-1}^*<1 $ such that for  $j=1,\cdots,K^{*}$, $t_j^{*}=[n\tau_j^*]$ (where $[x]$ is the floor of $x$). The vector $\underline{\tau}^*=\big(\tau_1^*,\ldots,\tau_{K^*-1}^*\big) $ is called the vector of breaks.}\\
~\\
Even if the length of $T^*_j$ has asymptotically the same order than
$n$, the dependences with respect to $n$ of
 $t^*_j$, $t_k$, $T^*_j$ and $T_k$ are omitted for notational
 convenience. \\
Finally we make a technical non classical assumption. Using the convention: if $A_i(M_\theta,\Theta)$
 holds then $\alpha_\ell^{(i)}( h_\theta,\Theta) = 0$ and if $A_i(h_\theta,\Theta)$
 holds then $\alpha_\ell^{(i)}( M_\theta,\Theta) = 0$, define:\\
~\\
{\bf Hypothesis H$_i$ ($i=0,1,2$)}: {\em For $0\le p\le i$, the assumptions $A_p(f_\theta,\Theta)$,
$A_p(M_\theta,\Theta)$ (or $ A_p(h_\theta,\Theta)$) hold and for all
$j=1,\ldots,K^*$ there exists $r\ge1$ such that $\theta_j^* \in
\Theta(r)$. Denoting
$$
c^*=  \underset{j=1,\cdots,K^*} {\mbox{min}}\big(-\log
(\beta^{(0)}(\theta^*_j))/8 \big)\wedge \underset{j=1,\cdots,K^*}
{\mbox{min}}\big(-\log (\tilde\beta^{(0)}(\theta^*_j))/8 \big)
$$
the  sequence $(v_n)_{n\in \N}$ used in
 \eqref{contrast}  satisfies   for all  $j=1,\cdots K^*$:
\begin{multline}
  \sum_{k\ge 1}\Big(\dfrac{v_k}{k}\Big)^{r/4\wedge1} \Big(\hspace{-5mm}\sum\limits_{\ell\ge kc^*/\log(k)}
 \beta_\ell^{(0)}(\theta^*_j)\Big)^{r/4} \hspace{-3mm}\wedge  \sum_{k\ge 1}\Big(\dfrac{v_k}{k} \Big)^{r/4\wedge1} \Big(\hspace{-5mm}\sum\limits_{\ell\ge kc^*/\log(k)}
 \tilde\beta_\ell^{(0)}(\theta^*_j)\Big)^{r/4}\hspace{-5mm} < \hspace{-1mm}\infty~~~ \mbox{and}\\
 \label{condv}
 \sum_{k\ge 1}\Big(\dfrac{v_k}{k}\Big)^{r/4\wedge1} \Big( \sum\limits_{\ell\ge
 k  /2} \big(\alpha_\ell^{(p)}( f_\theta,{\Theta(r)})+ \alpha_\ell^{(p)}( M_\theta,{\Theta(r)}) + \alpha_\ell^{(p)}( h_\theta,{\Theta(r)})\big)\Big)^{r/4}\hspace{-5mm}< \hspace{-1mm}\infty.
\end{multline}}
\noindent The assumption H$_i$ is interesting as it links the
decrease rate of the Lipschitz coefficients and the penalization
term of \eqref{contrast}. The classical BIC penalization and the one
coming from the MDL approach (see \cite{David2008}) correspond to a
sequence $v_n\propto n/\log(n)$. This choice is possible if the
Lipschitz coefficients decrease exponentially fast, which hold for
all models {\cal M}$(f_\theta,M_\theta)$ with finite order (see
Remark below). However, if the decrease of the Lipschitz
coefficients is slower, our method can exclude such a choice and an heavier term $\beta_n=n/v_n>>\log(n)$ in the penalization has to be chosen.
\begin{rem}\label{ratevn}
Conditions \eqref{condv} satisfied by $(v_n)_n$ are deduced from a
result of Kounias \cite{Kounias1969}. The conditions on $(v_n)_n$
are not too restrictive:
\begin{description}
    \item[(1)]  geometric case: if ~ $ \alpha_\ell^{(i)}( f_\theta,\Theta(r))+ \alpha_\ell^{(i)}( M_\theta,\Theta(r)) + \alpha_\ell^{(i)}( h_\theta,\Theta(r)) =O(a^\ell)$ with $0\leq a <1$,
    then any $(v_n)_n$ such as $v_n=o(n)$ can be chosen (for instance $v_n=n (\log n)^{-1}$).
    \item[(2)]  Riemanian case: if ~ $ \alpha_\ell^{(i)}( f_\theta,\Theta(r))+ \alpha_\ell^{(i)}( M_\theta,\Theta(r)) + \alpha_\ell^{(i)}( h_\theta,\Theta(r)) =O(\ell^{-\gamma})$ with $\gamma>1$,
    \begin{itemize}
    \item if $\gamma > 1+  (1 \vee 4r^{-1})$,
    then all sequence $(v_n)_n$ such as $v_n=o(n)$ can be chosen (for instance $v_n=n (\log n)^{-1}$).
    \item if $1 \vee 4r^{-1} < \gamma\le 1+(1 \vee 4r^{-1})$, then any $(v_n)_n$ such as $v_n=O(n^{\gamma-(1 \vee 4r^{-1})}(\log n)^{-\delta})$ with $\delta >1 \vee 4r^{-1}$ can be chosen.
    \end{itemize}
\end{description}
\end{rem}
We are now ready to prove the consistency of the penalized QMLE:
\begin{Theo}\label{theo4}
 Assume that the hypothesis $D({\Theta(r)})$,
$Id({\Theta(r)})$, B, C and H$_0$ are satisfied with $r\ge 2$ and
$v_n\to \infty$.  If $K_{\max} \ge K^\ast$  then:
\begin{eqnarray}\label{consist}
(\widehat{K}_n, \widehat{\underline{\tau}}_n,
\widehat{\underline{\theta}}_n) \limiteproban
(K^*,\underline{\tau}^*,\underline{\theta}^*).
\end{eqnarray}
\end{Theo}

\begin{rem}
If $K^*$ is known, we can relax the assumptions for the consistency
by taking $v_n=1$ for all $n$ as the penalization term in
\eqref{contrast} does not matter. If $K^*$ is unknown then a reasonable choice in any geometric or Riemanian cases is $v_n\propto\log n$ (therefore $\beta_n\propto n (\log n)^{-1}$), see Remark \ref{ratevn}. 
\end{rem}

\subsection{Rate of convergence of the estimators}
To state a rate of convergence of the estimators
$\widehat{\underline{\tau}}_n$ and $\widehat{\underline{\theta}}_n$,
we need to work under stronger moment and regularity assumptions.
\begin{Theo}\label{theo3}
Assume that the hypothesis $D({\Theta(r)})$, $Id({\Theta(r)})$, B, C
and H$_2$ are satisfied with $r\ge 4$ and $v_n\to \infty$.  If
$K_{\max} \ge K^\ast$  then  the sequence
    $( \| \widehat{\underline{t}}_n - \underline{t}^* \|_m )_{n>1}$ is uniformly tight in probability, {\em i.e.}
 \begin{equation}\label{tight} ~~ \mathop {\lim }\limits_{\delta \, \to \, \infty } ~ \mathop {\lim
  }\limits_{n \, \to  \infty } \mathbb{P}( \| \widehat{\underline{t}}_n - \underline{t}^* \|_m > \delta ) = 0.  \end{equation}
\end{Theo}
This theorem induces that $w_n^{-1} \, \| \widehat{\underline{t}}_n -\underline{t}^*\|_m \overset{P}\rightarrow 0  $ for any sequence $(w_n)_n$ such as $w_n \to \infty$ and therefore $ \| \widehat{\underline{t}}_n -\underline{t}^*\|_m=o_P(w_n)$: the convergence rate is arbitrary close to $O_P(1)$. This is the same convergence rate as in the case where $(X_t)_t$ is a sequence of independent r.v. (see for instance \cite{Bai1998}). Such convergence rate was already reached for mixing processes in \cite{LavielleMoulines2000}. \\
~\\
Let us turn now the convergence rate of the estimator of parameters
$\theta_j^*$. By convention if $\widehat K_n< K^*$, set $\widehat T_j=\widehat T_{\widehat K_n}$ for $j\in \{\widehat K_n,\ldots,K^*\}$. Then,
\begin{Theo}\label{TLC}
Assume that the hypothesis $D({\Theta(r)})$, $Id({\Theta(r)})$, B, C
and H$_2$ are satisfied with $r\ge 4$ and $\sqrt n=O(v_n)$. Then if
$\theta^*_j\in\overset{\circ}{\Theta}(r)$ for all $j=1,\cdots,K^*$,
we have
\begin{eqnarray}\label{tlc}
  \sqrt{n^{*}_j} \, \big( \widehat{\theta}_n(\widehat{T}_j) - \theta^*_j \big) \overset{\mathcal{ D}}{\underset{n\to \infty} \longrightarrow}
   \mathcal{N}_d \big(0, F(\theta^*_j)^{-1}G(\theta^*_j)F(\theta^*_j)^{-1} \big),
\end{eqnarray}
where, using $q_{0,j}$ defined in \eqref{def_qj}, the matrix $F$ and
$G$ are such as
  \begin{equation}\label{eq_AN}
  (F(\theta^*_j))_{k,l} = \E \Big( \dfrac{\partial^2 q_{0,j}(\theta^*_j)}{\partial \theta_k \partial \theta_l}
  \Big) ~ \text{and} ~~ (G(\theta^*_j))_{k,l} = \E \Big( \dfrac{\partial q_{0,j}(\theta^*_j)}{\partial \theta_k} \dfrac{\partial
  q_{0,j}(\theta^*_j)}{\partial \theta_l}\Big).
  \end{equation}
\end{Theo}
\begin{rem}\label{choicevn}
In Theorem \ref{TLC}, a condition on the rate of convergence of
$v_n$ is added. The optimal choice for the penalization term corresponds to $v_n\propto \sqrt n$ as it corresponds to the most general problem \eqref{problem}, see Remark \ref{ratevn}. However, by assumption H$_2$ it excludes models with finite
moments $r\ge 4$ satisfying:  $\ell^{-\gamma} =O( \alpha_\ell^{(i)}( f_\theta,\Theta(r))+
\alpha_\ell^{(i)}( M_\theta,\Theta(r)) + \alpha_\ell^{(i)}(
h_\theta,\Theta(r)))$ with $1<\gamma\le 3/2$ for
some $i=0,1,2$. For these models the consistency and the rate of convergence of order $n$ for $\widehat {\underline\tau}_n$ hold but we do not get any rate of convergence for $\widehat {\underline \theta}_n$.
\end{rem}

\section{Some examples} \label{Examples}

\subsection{AR($\infty$) models}
Consider $AR(\infty)$ with $K^*-1$ breaks defined by the equation:
$$ X_t = \sum_{k\geq 1} \phi_k(\theta^*_j)X_{t-k} + \xi_t  \ ,\ \ t^*_{j-1}< t \leq t^*_j \ , \ j=1,\cdots,K^*. $$
It correponds to the problem \eqref{problem} with models ${\cal
M}_{T^*_i}(f_\theta,M_\theta)$ where  $f_\theta(x_1,\cdots)= \sum_{k\geq 1} \phi_k(\theta)x_k $ and $ M_\theta \equiv 1 $.
Assume that $\Theta$ is a compact set such that $\sum_{k\geq 1}
\|\phi_k(\theta)\|_\Theta <1$. Thus  $\Theta(r)=\Theta$ for any
$r\ge1$ satisfying $\E|\xi_0|^r<\infty$. Then Assumptions
D($\Theta$) and $A_0(f_\theta,\Theta)$ hold automatically with
$\underline{h}=1$ and $\alpha^{(0)}_k(f_\theta,\Theta(r))=
\|\phi_k(\theta)\|_\Theta$. Then,
\begin{itemize}
\item Assume that Id($\Theta$) holds and that there exists $r\geq 2$ such that $\E|\xi_0|^r<\infty$. If there exists $\gamma >1 \vee 4r^{-1}$ such that  $\|\phi_k(\theta)\|_\Theta = O(k^{-\gamma})$ for all $k \geq 1$,
then the penalization $v_n=\log n$ (or $\beta_n=n /\log n $) ensures the consistency of $(\widehat{K}_n, \widehat{\underline{\tau}}_n,
\widehat{\underline{\theta}}_n)$.
\item Moreover, if $r\ge 4$, $\gamma > 3/2 $ and $\phi_k$ twice
differentiable satisfying $ \|\phi_k'(\theta)\|_\Theta =
O(k^{-\gamma})$ and $ \|\phi_k''(\theta)\|_\Theta = O(k^{-\gamma})$,
then the penalization $v_n=\beta_n=\sqrt n$ ensures the convergence (\ref{tight}) of  $\widehat{\underline{t}}_n$ and the CLT (\ref{tlc}) satisfied by $\widehat{\underline \theta}_n(\widehat T_j)$ for all $j$.
\end{itemize}
Note that this
problem of change detection was considered by Davis {\it et
al.} in \cite{David1995} under moments
of order greater than $4$ is  required. In Davis {\it
et al.} \cite{David2008}, the same problem for another  break model for AR processes is
studied. However, in both these papers, the process is supposed to be independent from one
block to another and stationary on each block.

\subsection{ARCH($\infty$) models}
 Consider  an $ARCH(\infty)$ model  with  $K^*-1$ breaks defined by:
$$
 X_t=\Big ( \psi_0(\theta^*_j)+\sum\limits_{k=1}^{\infty}\psi_k(\theta^*_j)X^2_{t-k} \Big )^{1/2}\,
 \xi_t \ ,\ \ t^*_{j-1}< t \leq t^*_j \ , \ j=1,\cdots,K^*,
$$
 where for any $\theta \in \Theta$, $\psi_0(\theta)>0$ and $(\psi_k(\theta))_{k\ge 1}$ is a sequence of positive real number and
 $\E(\xi_0^2)=1$. Note that $ h_{\theta} ((x_k)_{k\in \mathbb{N}} )=
 \psi_0(\theta)+\sum_{k=1}^{\infty}\psi_k(\theta)x^2_{k} $ \ and $f_\theta= 0.$  Assume that $\Theta$ is a compact set such that $\sum_{k\geq 1}
 \|\psi_k(\theta)\|_\Theta <1$, then $\Theta(2)=\Theta$. Assume   that $ \inf_{ \theta\in \Theta }\psi_0(\theta)>0$ which ensures that $D(\Theta)$
 and $Id(\Theta)$ hold.
\begin{itemize}
\item If there exists $\gamma >2$ such that $\|\psi_k(\theta)\|_\Theta = O(k^{-\gamma})$ for all $k \geq 1$, then the penalization $v_n=\log n$ (or $\beta_n= n/\log n$) leads to the consistency of $(\widehat{K}_n, \widehat{\underline{\tau}}_n,
\widehat{\underline{\theta}}_n)$ when $\theta^*_j \in \Theta$ for all $j$.
\item Moreover, if $r\ge 4$ and $\psi_k$ is twice differentiable satisfying
 $\|\psi_k'(\theta)\|_\Theta = O(k^{-\gamma})$ and $
 \|\psi_k''(\theta)\|_\Theta = O(k^{-\gamma})$ with $\gamma > 3/2$, if $\Theta(4)$ is a compact such that $\theta^*_j \in\overset{\circ}{\Theta}(4)$ for all $j$,
then the penalization $v_n=\beta_n=\sqrt n $ as in Remark \ref{ratevn}  ensures the convergence (\ref{tight}) of  $\widehat{\underline{t}}_n$ and the CLT (\ref{tlc}) satisfied by $\widehat{\underline \theta}_n(\widehat T_j)$ for all $j$.
\end{itemize}
This problem of break detection was already studied by Kokoszka and Leipus in
 \cite{Leipus2000} but they obtained the consistency of their procedure under stronger assumptions.
\begin{Ex}
Let us detail the
GARCH$(p,q)$ model  with  $K^*-1$ breaks defined by:
$$ X_t= \sigma_t \, \xi_t \ ,\ \sigma_t^2 = a_{0,j}^*+ \sum^{q}_{k=1} a_{k,j}^*X^2_{t-k} +
  \sum^{p}_{k=1} b_{k,j}^*\sigma^2_{t-k}  \,  \ \  t^*_{j-1}< t \leq t^*_j \ , \ j=1,\cdots,K^*
 $$
 with $\E(\xi_0^2)=1$. Assume that for any $\theta=(a_0,\ldots,a_q,b_1,\ldots,b_p) \in \Theta $ then  $a_k\geq0$, $b_k\geq0$ and $ \sum_{k=1}^{p}b_k <1  $.
 Then, there exists (see Nelson and Cao \cite{Nelson1992}) a nonnegative
 sequence $(\psi_k(\theta))_k $ such that $\sigma^2_t= \psi_0(\theta)+ \sum_{k\geq 1}\psi_k(\theta)X^2_{t-k}.$ Remark that this sequence is twice differentiable with respect to $\theta$ and that its derivatives are exponentially decreasing.
Moreover for any $\theta \in \Theta $ it holds
   $\sum_{k\geq 1}\psi_k(\theta) \leq \big( \sum_{k=1}^{q}a_k\big)/\big(1- \sum_{k=1}^{p}b_k\big)$ and one can consider:
$$
\Theta(r)= \left\{ \theta \in \Theta ~, ~  (\E | \xi_0|^r)^{2/r} \sum_{k=1}^{q}a_k +  \sum_{k=1}^{p}b_k <1 \right\}.
$$
Then if  $\sum_{k=1}^{q}a_{k,j}^* +  \sum_{k=1}^{p}b_{k,j}^* <1$ for all $j$ (case $r\geq 2$), our estimation procedure associated with a penalization term $\beta_n K$ for any $1<<\beta_n<<n$ is consistent. Moreover, if $(\E | \xi_0|^4)^{1/2}\sum_{k=1}^{q}a_{k,j}^* +  \sum_{k=1}^{p}b_{k,j}^* <1$ for all $j$, then our procedure with a penalization $1<<\beta_n=0(\sqrt n)$ allows the same rates of convergence than in the case where $(X_t)$ are independent r.v. For example, a penalization $\beta_n \propto \log n$ as in \cite{David2008} can be chosen in this case.
 \end{Ex}

\subsection{Estimates breaks in TARCH$(\infty$) model}
Consider a TARCH($\infty$) model with breaks defined by:
$$ X_t= \sigma_t \, \xi_t \ ,\ \sigma_t = b_0(\theta^*_j)+ \sum_{k\geq1}
 \Big(b^+_k(\theta^*_j)\text{max}(X_{t-k},0)- b^-_k(\theta^*_j)\text{min}(X_{t-k},0)  \Big),
$$
for any $ t^*_{j-1}< t \leq t^*_j, \ j=1,\cdots,K^*$ and where $\sum_{k\geq1}\text{max}( \| b_k^+(\theta) \|_\Theta,
\| b_k^-(\theta) \|_\Theta)<\infty$. Then $f_\theta=0$ and
($A_0(M_\theta,\Theta)$) holds with
$\alpha^{(0)}_k(M_\theta,\Theta)= \text{max}( \| b_k^+(\theta)
\|_\Theta,  \| b_k^-(\theta) \|_\Theta) $.
\begin{itemize}
\item Assume  $ \| \xi_0
\|_r\sum_{k\geq1}\text{max}( \| b_k^+(\theta) \|_\Theta,  \|
b_k^-(\theta)
 \|_\Theta)<1$ for $r\geq 2$. If there exists  $\gamma > 1 \vee 4r^{-1}$ such as $\text{max}( \| b_k^+(\theta) \|_\Theta,  \| b_k^-(\theta)
  \|_\Theta)=O(k^{-\gamma})$  for all $k \geq 1$, then a penalization $v_n=\log n$ (or $\beta_n=n /\log n$) leads to the consistency of $(\widehat{K}_n, \widehat{\underline{\tau}}_n,
\widehat{\underline{\theta}}_n)$ when $\theta^*_j\in\Theta(2)$ for all $j$.
\item Moreover, if $r\ge 4$ and $b^+_k$, $b^-_k$ are twice differentiable satisfying
 $\|\partial b^+_k(\theta)/\partial \theta \|_\Theta = O(k^{-\gamma})$ and $
 \|\partial^2 b^-_k(\theta)/\partial \theta^2\|_\Theta = O(k^{-\gamma})$ with $\gamma > 3/2$ (the same for $b^-_k$), then $v_n=\beta_n=\sqrt n$  ensures the convergence (\ref{tight}) of  $\widehat{\underline{t}}_n$ and the CLT (\ref{tlc}) satisfied by $\widehat{\underline \theta}_n(\widehat T_j)$ for all $j$ (with $\theta_j^* \in\overset{\circ}{\Theta}(4)$).
\end{itemize}


\section{Proofs of the main results}\label{proofs}
In the sequel $C$ denotes a positive constant whom value may differ
from one inequality to another.

\subsection{Proof of Proposition \ref{prop1}}
 (i) It is clear that $\{X_t,\, t\leq t^{*}_1\}$ exists and is causal, stationary with finite moments of order $r$ (see \cite{Bardet2009}). Therefore, $X$ is defined by induction as follows:
\begin{equation}\label{eq2}
 X_t:=M_{\theta^*_j}(X_{t-1},X_{t-2},\cdots) \xi_t+ f_{\theta^*_j}(X_{t-1},X_{t-2},\cdots), \;\forall t\in T^{*}_j ; \;
 j=2,\cdots K^{*}.
\end{equation}
Thus, $X_t$ is independent of $(\xi_j)_{j>t}$ and it suffices to
prove (ii) which immediately leads  existence of moments.

(ii) Let us first consider the general case when
$A_0(f_\theta,\{\theta\})$ and $A_0(M_\theta,\{\theta\})$ hold with
$\beta^{(0)}(\theta)<1$. As in \cite{Doukhand}  we remark that
$$
\|X_t\|_r\le\frac{\|Z_{t^\ast_j,1}\|_r}{1-\beta^{(0)}(\theta^*_1)}
$$
 for $t\le t_1^\ast$, with $Z_{t,j}:=M_{\theta^*_j}(0,0,\cdots) \xi_t+ f_{\theta^*_j}(0,0,\cdots)$ for all $j=1,\ldots,K^*$.  Assume that there exists $C_{r,t}>0$ such that $C_{r,t}=\sup_{i<t}\|X_i\|_r$  and let $t\in T^\ast_j$,
then
$$
     |X_t-Z_{t,j}| \leq |M_{\theta^*_j}(X_{t-1},\cdots)-M_{\theta^*_j}(0,0,\cdots)| |
     \xi_t|+|f_{\theta^*_j}(X_{t-1},\cdots)-f_{\theta^*_j}(0,0,\cdots)|.
     $$
      We obtain for all $ t $, by independence of  $(\xi_j)_{j>t}$  and $ X_t $:
     $$
     \|X_t-Z_t\|_r \leq \|M_{\theta^*_j}(X_{t-1},\cdots)-M_{\theta^*_j}(0,0,\cdots)\|_r   \|
     \xi_t\|_r+\|f_{\theta^*_j}(X_{t-1},\cdots)-f_{\theta^*_j}(0,0,\cdots)\|_r.
     $$
     Then, we have:
\begin{align*}
     \|M_{\theta^*_j}(X_{t-1},\cdots)-M_{\theta^*_j}(0,0,\cdots)\|_r
    & \leq  \sum\limits_{i=1}^{\infty} \alpha^{(0)}_i(M_{\theta^*_j},\theta^*_j)\|X_{t-i}\|_r \leq C_{r,t}\sum\limits_{i=1}^{\infty} \alpha^{(0)}_i(M_{\theta^*_j},\theta^*_j),\\
        \|f_{\theta^*_j}(X_{t-1},\cdots)-f_{\theta^*_j}(0,0,\cdots)\|_r
    & \leq  \sum\limits_{i=1}^{\infty} \alpha^{(0)}_i(f_{\theta^*_j},\theta^*_j)\|X_{t-i}\|_r \leq C_{r,t}\sum\limits_{i=1}^{\infty} \alpha^{(0)}_i(f_{\theta^*_j},\theta^*_j).
\end{align*}
We deduce that
$$
\|X_t\|_r\le
\|Z_{t,j}\|_r+C_{r,t}\left(\sum\limits_{i=1}^{\infty}\alpha^{(0)}_i(f_{\theta^*_j},\{\theta^*_j\})+(\E\|\xi_0\|^r)^{1/r}\sum\limits_{i=1}^{\infty}
\alpha^{(0)}_i(M_{\theta^*_j},\{\theta^*_j\})\right).
$$
Thus, $\|X_t\|^r<\infty$, $C_{r,t+1}<\infty$ and $ \|X_t\|_r\le
\|Z_{t,j}\|_r+C_{r,t+1}\beta^{(0)}(\theta^*_j) $ since $C_{r,t}\le
C_{r,t+1}$. Similarly for any $i<t$, we have $C_{r,i}\le C_{r,t+1}$
 and $
\|X_i\|_r\le \max_{1\le j\le
K^\ast}\left\{\|Z_{t,j}\|_r+C_{r,t+1}\beta^{(0)}(\theta^*_j)\right\}$.
Thus, by definition of $C_{r,t+1}=\sup_{i\le t}\|X_t\|_r$ we obtain
$$
C_{r,t+1}\le \max_{1\le j\le
K^\ast}\left\{\|Z_{t,j}\|_r+C_{r,t+1}\beta^{(0)}(\theta^*_j)\right\},
$$
and the Proposition is established.\\

In the ARCH-type case when $f_\theta=0$ and
$A_0(h_\theta,\{\theta\})$ holds with $\tilde\beta^{(0)}(\theta)<1$,
we follow the same reasoning than previously starting from the
inequality
$$
\|X^2_t-(M_{\theta^*_j}(0,0,\cdots) \xi_t)^2\|_{r/2} \leq
\|h_{\theta^*_j}(X_{t-1},\cdots)-h_{\theta^*_j}(0,0,\cdots)\|_{r/2}
\|\xi_t^2\|_{r/2}.$$ Finally we obtain the desired result  with
$$
C= \max_{1\le j \le K^\ast}\frac{\|M_{\theta^*_j}(0,0,\cdots)\
\xi_0+ f_{\theta^*_j}(0,0,\cdots)\|_r
}{1-\beta^{(0)}(\theta^*_j)}\wedge \max_{1\le j \le
K^\ast}\frac{\|M_{\theta^*_j}(0,0,\cdots)\ \xi_0\|_r
}{\sqrt{1-\tilde\beta^{(0)}(\theta^*_j)}}.~~~~~~ \blacksquare
$$

\subsection{Some preliminary result} The following technical  lemma is useful in the sequel:
\begin{lem} \label{lem1}
Suppose that $\theta_j^*\in\Theta(r)$ for $ j=1,\ldots,K^*$ with
$r\ge 2$ and under the assumptions $A_0(f_\theta,\Theta)$,
$A_0(M_\theta,\Theta)$ (or $A_0(h_\theta,\Theta)$) and
$D(\Theta(r))$, then there
 exists $C>0$ such that
$$
\mbox{for all $t\in \Z$,}\quad
 \mathbb{E}\Big(\sup_{\theta\in\Theta(r)}\big|q_t(\theta)\big|\Big)\le C.
$$
\end{lem}
{\bf Proof} Using the inequality $(a+b)^2\leq 2(a^2+b^2)$, we have for
all $ t\in \Z$:
\begin{eqnarray*}
  \big\|f^t_\theta \big \|^2_{\Theta(r)} &\le &2 \Big(\big\| f^t_\theta-f_\theta(0,\ldots)\big\|_{\Theta(r)}^2+\big\|f_\theta(0,\ldots)\big\|_{\Theta(r)}^2 \Big)  \\
   &\le& 2 \Big( \Big(\sum_{i\ge 1}\alpha_i^{(0)}(f_\theta,{\Theta(r)})\Big) \cdot \sum_{i\ge 1}\alpha_i^{(0)}(f_\theta,{\Theta(r)})|X_{t-i}|^2
+\big\|f_\theta(0,\ldots)\big\|_{\Theta(r)}^2 \Big),
\end{eqnarray*}
therefore
$$
 \mathbb{E}\big\|f^t_\theta \big\|^2_{\Theta(r)} \leq 2 \Big ( C\Big(\sum_{i\ge 1}\alpha_i^{(0)}(f_\theta,{\Theta(r)})\Big)^2
 +\big\|f_\theta(0,\ldots)\big\|_{\Theta(r)}^2 \Big).
$$
Thus $ \mathbb{E}\big\|f^t_\theta\|^2_{\Theta(r)} \le C $  for all $t\in \mathbb{Z}$ and similarly  $\mathbb{E}(\|h^t_\theta\|_{\Theta(r)})=\mathbb{E}(\|M^t_\theta\|^2_{\Theta(r)}) \le C_M$. Yet, under  assumption
 $(D({\Theta(r)}))$, we have: $|q_t(\theta)|\le \dfrac{1}{\underline{h}}|X_t- f_\theta^t|^2+|\log(h_\theta^t)| $
 and using inequality $\log x\le x-1$ for all $x>0$, it follows:
$$
 |\log(h_\theta^t)|=\Big|\log(\underline{h})+\log(\frac{h_\theta^t}{\underline{h}})\Big|
 \le  1+|\log(\underline{h})|+\frac{1}{\underline{h}}h_\theta^t.
$$
Finally, we have for all $t\in\Z$:
\begin{eqnarray*}
  \mathbb{E}\big(\sup\limits_{\theta\in {\Theta(r)} }|q_t(\theta)|\big) &\le& 1+|\log \underline{h}|
  +\frac{1}{\underline{h}}\big(\mathbb{E}\|h^t_\theta \|_{{\Theta(r)}}+ 2\mathbb{E}|X_t|^2+ 2\mathbb{E}\|f^t_\theta \|_{{\Theta(r)}}^2\big) \le C.~~~~~~~ \blacksquare
\end{eqnarray*}

\subsection{Comparison with  stationary solutions}
In the following, we assume that $\theta^*_j\in \Theta(r)$ for all
$j=1,\ldots, K^*$ with $r\ge 1$. It comes from
 \cite{Bardet2009} that the equation
$$
X_{t,j}=M_{\theta_j^*}\big ((X_{t-k,j})_{k\in \N^*}\big )\cdot\xi_t+ f_{\theta_j^*}\big ((X_{t-k,j})_{k\in \N^*}\big ) ~~ \text{for all} ~ t \in \Z
$$
 has  $r$ order stationary solution $\big(X_{t,j}\big)_{t\in \Z}$ for any $ j=1,\ldots, K^*$. Then
\begin{lem} \label{lem2}
Assume that  the assumptions $A_0(f_\theta,\Theta)$,
$A_0(M_\theta,\Theta)$ (or $A_0(h_\theta,\Theta)$) hold and that
$\theta_j^*\in\Theta(r)$ for $ j=1,\ldots,K^*$  for $r\ge 2$. Then:
\begin{enumerate}
\item $X_t=X_{t,1}$ for all $t\le t_1^*$;
\item There exists  $C>0$ such that for any $j\in \{2,\cdots,K^*\}$, for all  $t\in T^*_j$,
\begin{eqnarray*}\label{ineg}
 \|X_t-X_{t,j}\|_{r}&\le& C\Big(\inf_{1\le p\le t-t^*_{j-1}}\Big\{\beta^{(0)}(\theta_j^\ast)^{(t-t^*_{j-1})/p}+
 \sum_{i\ge p}\beta_i^{(0)}(\theta_j^\ast)\Big\}\Big)\\
 \label{inegbis}  \|X_t^2-X_{t,j}^2\|_{r/2}&\le& C
 \Big( \inf_{1\le p\le t-t^*_{j-1}}\Big\{\tilde\beta^{(0)}(\theta_j^\ast)^{(t-t^*_{j-1})/p}+
 \sum_{i\ge p}\tilde\beta_i^{(0)}(\theta_j^\ast)\Big\}\Big).
\end{eqnarray*}
\end{enumerate}
\end{lem}
\begin{dem}
 1. It is obvious from the definition of $X$.\\ 
2. Let $j\in \{ 2,\cdots,K^*\}$, we proceed by induction on $t\in T^*_j$.\\
First consider the general case where $A_0(f_\theta,\{\theta\})$ and
$A_0(M_\theta,\{\theta\})$ hold with $\beta^{(0)}(\theta)<1$. By
Proposition \ref{prop1}, there exists $C_r\geq 0$ such that $
\|X_t^2-X_{t,j}^2\|_{r/2}\le \|X_t\|_r+ \|X_{t,j} \|_r \le C+
\max_{1\le j \le K^\ast}\|X_{0,j} \|_r\le C_r$ for all $
j=1,\ldots,K^*$ and $t\in\Z$. For $1\le p\le t-t^*_{j-1}$ let
$u_\ell:=\sup_{t^*_{j-1}+\ell p\le i\le t^*_j}\|X_i-X_{i,j}\|_r$. Then
$\|X_t-X_{t,j} \|_r\le u_{[(t-t^*_{j-1})/p]}$ and for any $t\le i\le
t^*_{j}$:
\begin{align*}
\|X_{i}-X_{i,j}\|_r &\le \sum\limits_{k\ge 1} \beta^{(0)}_{k}(\theta_j^*)\|X_{i-k}-X_{i-k,j}\|_r\\
 &\le \sum_{k=1}^p \beta^{(0)}_{k}(\theta_j^*)\|X_{i-k}-X_{i-k,j}\|_r+C_r  \sum_{k>p}\beta^{(0)}_{k}(\theta_j^*)\\
 &\le \beta^{(0)}(\theta_j^*)u_{[(t-t^*_{j-1})/p]-1}+C_r \sum_{k>p}\beta^{(0)}_{k}(\theta_j^*).
\end{align*}
Similarly, it is easy to show that for all $1\le \ell\le
[(t-t^*_{j-1})/p]$ we have
$$
 u_\ell\le \beta^{(0)}(\theta_j^*)u_{\ell-1}+C_r \sum_{k>p}\beta^{(0)}_{k}(\theta_j^*).
$$
Denote $a=\beta^{(0)}(\theta_j^*)<1,  \ b= C_r
\sum_{k>p}\beta^{(0)}_{k}(\theta_j^*)$  such that $u_\ell \leq a
 u_{\ell-1} + b$. Considering $  w_0=u_0 \ \text{and} \ w_l=aw_{\ell-1}+b$, then $w_\ell=a^\ell w_0+b(1-a^{\ell-1})/(1-a)
\leq a^\ell w_0 + {b}/({1-a})$.  Since $u_0 \leq C_r$ by definition
and  $u_\ell \leq w_{\ell}$ for any $\ell$, we have:
\begin{eqnarray*}
u_\ell  &\leq &a^\ell u_0 + \dfrac{b}{1-a} \leq \big(  \beta^{(0)}(\theta_j^*) \big) ^\ell C_r + \dfrac{C_r}{1-\beta^{(0)}(\theta_j^*)}\sum_{k>p}\beta^{(0)}_{k}(\theta_j^*)\\
&\leq &\dfrac{C_r}{1-\beta^{(0)}(\theta_j^*)} \Big(
\beta^{(0)}(\theta_j^*) \big) ^\ell +
\sum_{k>p}\beta^{(0)}_{k}(\theta_j^*) \Big).
\end{eqnarray*}
Thus for all $1\le p\le t-t^*_{j-1}$
$$
 \|X_t^2-X_{t,j}^2\|_{r/2} \le C_r\|X_t-X_{t,j}\|_r  \le C_r u_{[(t-t^*_{j-1})/p]}\le C \big(\beta^{(0)}(\theta^*_j)^{(t-t^*_{j-1})/p}+\sum_{i\ge
p}\beta_i^{(0)}(\theta^*_j)\big)
$$
and Lemma \ref{lem2} is proved.\\
~\\
In the ARCH-type case when $f_\theta=0$ and
$A_0(h_\theta,\{\theta\})$ holds with $\tilde\beta^{(0)}(\theta)<1$,
we follow the same reasoning than previously starting from the
inequality
$$
\|X_{i}^2-X_{i,j}^2\|_{r/2}  \le \sum\limits_{k\ge 1}
\tilde\beta^{(0)}_{k}(\theta_j^*)\|X_{i-k}^2-X_{i-k,j}^2\|_{r/2}.
$$
For all $ j=1,\ldots,K^*$ and $t\in \mathbb{Z}$, by Proposition
\ref{prop1}, $\|X_{i}^2-X_{i,j}^2\|_{r/2}\le  C_r^2$ and therefore
$$
\tilde u_\ell\le \tilde\beta^{(0)}(\theta_j^*)\tilde
u_{\ell-1}+C_r^2 \sum_{k>p}\tilde\beta^{(0)}_{k}(\theta_j^*)
$$
for $\tilde u_\ell=\sup_{t^*_{j-1}+\ell p\le i\le
t^*_j}\|X_i^2-X_{i,j}^2\|_{r/2}$ and Lemma \ref{lem2} is proved.
\end{dem}
\vspace{-1cm}
\subsection{The asymptotic behavior of the likelihood}
For the process $\big( X_{t,j}\big )_{t\in T^*_j,\,j=1,\ldots,
K^*}$, for any $j\in \{1,\cdots,K^*\}$ and $s\in T^*_j$ denote:
\begin{equation}\label{def_qj}
 q_{s,j}(\theta):=\frac{\left(X_{s, j}-f_\theta^{s,j}\right)^2}{h^{s,j}_{\theta}}+\log \left(h^{s,j}_\theta\right)
\end{equation}
with  $f_\theta^{s,j}:=f_\theta(X_{s-1,j},X_{s-2,j},\ldots)$,
$h^{s,j}_{\theta}:=(M_\theta^{s,j})^2$
 where $M_\theta^{s,j}:=M_\theta(X_{s-1,j},X_{s-2,j},\ldots)$. For any $ T\subset T_j^*$, denote
$$
L_{n,j}(T,\theta):=-\frac{1}{2}\sum\limits_{s\in T}q_{s,j}(\theta)
$$
 the likelihood of the  $j^{\textrm{th}}$ stationary model computed on  $T$.
\begin{lem}\label{lem_conv_L}
Assume that the hypothesis $D({\Theta(r)})$ holds.
\begin{enumerate}
\item If the assumption H$_0$ with $r\ge 2$ holds then for all $ j=1,\ldots,K^*$:
 $$
 \frac{v_{n^\ast_j}}{n_j^*}\left\|L_n\big(T^*_j,\theta\big)-L_{n,j}\big(T^*_j,
 \theta\big)\right\|_{\Theta(r)}\overset{\texttt{a.s.}}{\underset{n\to \infty}\longrightarrow}0.
$$
 \item  For $i=1,2$, if the assumption H$_i$ with $r\ge 4$ holds then for all $ j=1,\ldots,K^*$:$$
\frac{ v_{n^\ast_j}}{n_j^*}\Big \|\frac{\partial^i
L_n\big(T^*_j,\theta\big)}{\partial
 \theta^i}-\frac{\partial^i L_{n,j}\big(T^*_j,\theta\big)}{\partial
 \theta^i}\Big \|_{\Theta(r)}\overset{\texttt{a.s.}}{\underset{n\to \infty}\longrightarrow}0.
$$
\end{enumerate}

\end{lem}
\begin{dem}
 1-) For any $\theta\in \Theta(r)$,
$\Big|\frac{1}{n_j^*}L_n\big(T^*_j,\theta\big)-\frac{1}{n_j^*}L_{n,j}\big(T^*_j,
\theta\big)\Big | \le\dfrac{1}{n_j^*}\sum
\limits_{k=1}^{n_j^*}|q_{t^*_{j-1}+k}(\theta)-q_{t^*_{j-1}+k,j}(\theta)|$.
Then:
$$
  v_{n^*_j}\Big \|\frac{1}{n_j^*}L_n\big(T^*_j,\theta\big)-\frac{1}{n_j^*}L_{n,j}\big(T^*_j, \theta\big)\Big \|_{\Theta(r)} \le \frac{v_{n^*_j}}{n_j^*}\sum
\limits_{k=1}^{n_j^*}\|q_{t^*_{j-1}+k}(\theta)-q_{t^*_{j-1}+k,j}(\theta)\|_{\Theta(r)}.
$$
By Corollary 1 of Kounias \cite{Kounias1969},  with $r\le 4$ and no
loss of generality,  it is sufficient that
 $$
  \sum_{k\ge 1} \big(\frac{v_k }{k}\big)^{r/4}  \E\big(\|q_{t^*_{j-1}+k}(\theta)-q_{t^*_{j-1}+k,j}(\theta)\|^{r/4}_{\Theta(r)}\big)
  <\infty.
$$
For any $ \theta
\in \Theta(r)$, we have:
\begin{multline}\label{ineqdep} |q_s(\theta)-q_{s,j}(\theta)|\le \dfrac{1}{\underline{h}^2} |X_s- f_\theta^{s} |^2 |h^{s}_\theta-h^{s,j}_\theta | \\
+ \dfrac{1}{\underline{h}}
 \big( | X_s^2-X_{s,j}^2 | +  |f^{s}_\theta-f^{s,j}_\theta ||f^{s}_\theta+f^{s,j}_\theta+2X_s| + 2|f^{s,j}_\theta|| X_s-X_{s,j}  | + |h^{s}_\theta-h^{s,j}_\theta |\big).
 \end{multline}
First consider the general case with $A_0(f_\theta,\{\theta\})$ and
$A_0(M_\theta,\{\theta\})$ hold and $\beta^{(0)}(\theta)<1$:
\begin{multline*}
  \|q_s(\theta)-q_{s,j}(\theta)\|_{{\Theta(r)}} \le C \big( 1 + |X_{s,j}|+ |X_s|^2 +\|f_\theta^{s,j}\|_{\Theta(r)} +\|f_\theta^s\|^2_{\Theta(r)} \big)
 \\
  \times \big( | X_s-X_{s,j} |+ \|f^{s}_\theta-f^{s,j}_\theta\|_{\Theta(r)} +\|h^{s}_\theta-h^{s,j}_\theta\|_{\Theta(r)}\big),
\end{multline*}
and by Cauchy-Schwartz Inequality,
\begin{multline}
  (\E\|q_s(\theta)-q_{s,j}(\theta) \|_{{\Theta(r)}}^{r/4})^2\le C\E\big[\big( 1 + |X_{s,j}|+ |X_s|^2  +\|f_\theta^{s,j}\|_{\Theta(r)}
  +\|f_\theta^s\|^2_{\Theta(r)}\big)^{r/2}\big]\\
 \nonumber \times \E\big[\big( | X_s-X_{s,j} |+ \|f^{s}_\theta-f^{s,j}_\theta\|_{\Theta(r)} +\|h^{s}_\theta-h^{s,j}_\theta\|_{\Theta(r)}\big)^{r/2}\big].
\end{multline}
Using Proposition \eqref{prop1} and the argument of the proof of
Lemma \eqref{lem1} we claim that $\E|X_s|^r\le C$,
$\E\|f_\theta^s\|^r_{\Theta(r)}\le C$ and that
$\E\|f_\theta^{s,j}\|^r_{\Theta(r)}\le C$. Thus:
\begin{equation}\label{eq_q1}
 (\E\|q_s(\theta)-q_{s,j}(\theta) \|_{{\Theta(r)}}^{r/4})^2 \le C
 \big( \E | X_s-X_{s,j} |^{r/2} + \E\|f^{s}_\theta-f^{s,j}_\theta\|^{r/2}_{\Theta(r)} + \E\|h^{s}_\theta-h^{s,j}_\theta\|^{r/2}_{\Theta(r)}\big).
\end{equation}
Since $r/2\ge1$, we will use the $L^{r/2}$ norm. By Lemma
\ref{lem2}:
\begin{multline*} \|X_{s}-X_{s,j}\|_{r/2} \leq \|X_{s}-X_{s,j}\|_r  \le C\inf_{1\le p\le k}\big\{\beta^{(0)}(\theta_j^\ast)^{k/p}+\sum_{i\ge
 p}\beta_i^{(0)}(\theta_j^\ast)\big\}\\
 \le  C\inf_{1\le p\le k/2}\big\{\beta^{(0)}(\theta_j^\ast)^{k/(2p)}+\sum_{i\ge
 p}\beta_i^{(0)}(\theta_j^\ast)\big\}.
\end{multline*}
\vspace{-0.8cm} \begin{equation}\label{eq_x1}
 \Longrightarrow \quad \E|X_{s}-X_{s,j}|^{r/2} \leq C \Big( \inf_{1\le p\le k}\big\{\beta^{(0)}(\theta_j^\ast)^{k/p}+\sum_{i\ge
 p}\beta_i^{(0)}(\theta_j^\ast)\big\} \Big)^{r/2}.
\end{equation}
Moreover, as ($A_0(M_\theta,\Theta(r))$) holds, we have:
  \begin{equation}\label{eq_h1}
  \|\|h^s_\theta - h^{s,j}_\theta \|_{\Theta(r)}\|_{r/2}\le C \sum\limits_{i\geq 1}\alpha_i^{(0)}(M_\theta,{\Theta(r)})\|X_{s-i}-X_{s-i,j}\|_r.
\end{equation}
From \eqref{eq_h1} we obtain:
$$
  \|\|h^s_\theta - h^{s,j}_\theta \|_{\Theta(r)}\|_{r/2}\le C\Big( \sum\limits_{i=1}^{k/2-1}\alpha_i^{(0)}(M_\theta,{\Theta(r)})\|X_{s-i}-X_{s-i,j}\|_r
  +\sum\limits_{i\ge k/2}\alpha_i^{(0)}(M_\theta,{\Theta(r)})\|X_{s-i}-X_{s-i,j}\|_r\Big).
  $$
For all $s\ge t^*_{j-1}$ and $1 \leq i \leq k/2 - 1$, then $
s-i>t^*_{j-1}$, $ s-i> k/2$ and by Lemma \ref{lem2}:
\begin{eqnarray}
  \nonumber \|X_{s-i}-X_{s-i,j}\|_r &\leq&  C\inf_{1\le p\le k-i}\big\{\beta^{(0)}(\theta_j^\ast)^{(k-i)/p}+\sum_{i\ge p}\beta_i^{(0)}(\theta_j^\ast)\big\}  \\
 \nonumber  &\leq&  C\inf_{1\le p\le k/2}\big\{\beta^{(0)}(\theta_j^\ast)^{k/(2p)}+\sum_{i\ge p}\beta_i^{(0)}(\theta_j^\ast)\big\}
\end{eqnarray}
Thus, we can find $C>0$ not depending on $s$ such as:
\begin{equation}\label{eq_h2}
 \E\|h^s_\theta - h^{s,j}_\theta \|_{\Theta(r)}^{r/2}  \le C\Big( \inf_{1\le p\le k/2}\big\{\beta^{(0)}(\theta_j^\ast)^{k/(2p)}
 +\sum_{i\ge p}\beta_i^{(0)}(\theta_j^\ast)\big\}+\sum\limits_{i\ge k/2}\alpha_i^{(0)}(M_\theta,{\Theta(r)})\Big)^{r/2} .
\end{equation}
Similarly, we obtain:
\begin{equation}\label{eq_f1}
 \E\|f^s_\theta - f^{s,j}_\theta \|_{\Theta(r)}^{r/2}  \le C\Big( \inf_{1\le p\le k/2}\big\{\beta^{(0)}(\theta_j^\ast)^{k/(2p)}
 +\sum_{i\ge p}\beta_i^{(0)}(\theta_j^\ast)\big\}+\sum\limits_{i\ge k/2}\alpha_i^{(0)}(f_\theta,{\Theta(r)})\Big)^{r/2} .
\end{equation}
Relations \eqref{eq_q1}, \eqref{eq_x1}, \eqref{eq_h2} et
\eqref{eq_f1} give (the same inequality holds with $h_\theta$
replaced by $M_\theta$):
\begin{multline}\label{eqfin}
 \E\ \|q_s(\theta)-q_{s,j}(\theta)\|_{{\Theta(r)}}^{r/4} \le C\Big[  \Big( \inf_{1\le p\le k/2}\big\{\beta^{(0)}(\theta_j^\ast)^{k/(2p)}
  +\sum_{i\ge p}\beta_i^{(0)}(\theta_j^\ast)\big\} \Big)^{r/4}\\ +\Big(\sum\limits_{i\ge k/2}\alpha_i^{(0)}(f_\theta,{\Theta(r)})\Big)^{r/4}
+\Big(\sum\limits_{i\ge
k/2}\alpha_i^{(0)}(M_\theta,{\Theta(r)})\Big)^{r/4} \Big].
\end{multline}
By definition  $u_k = kc^*/\log(k)$ ($\leq k/2$ for large value of
$k$) satisfies the relation
$$\sum_{k\geq 1}\big( \dfrac{v_k}{k}\big)^{r/4}
 \big(\beta^{(0)}(\theta_j^\ast) \big)^{rk/8u_k}<\infty.
 $$
Choosing $p=u_k $ in \eqref{eqfin} we obtain:
\begin{multline*}
 \sum_{k\ge 1} \big(\frac{v_k }{k}\big)^{r/4}  \E\big(\|q_{t^*_{j-1}+k}(\theta)-q_{t^*_{j-1}+k,j}(\theta)\|^{r/4}_{{\Theta(r)}}\big)
  \leq  \sum_{k\geq 1}\big( \dfrac{v_k}{k}\big)^{r/4} \big(\beta^{(r)}(\theta_j^\ast) \big)^{rk/8u_k} \\
  \qquad +
  \sum_{k\geq 1}\big( \dfrac{v_k}{k}\big)^{r/4} \Big( \sum\limits_{i\ge u_k} \beta_i^{(0)}(\theta^*_j) \Big)^{r/4}
        +\sum_{k\geq 1}\big( \dfrac{v_k}{k}\big)^{r/4} \Big( \sum\limits_{i\ge
        k/2} \big(\alpha_i^{(0)}(f_\theta,{\Theta(r)})+\alpha_i^{(0)}(M_\theta,{\Theta(r)})\big)\Big)^{r/4}.
\end{multline*}
This bound is finite by assumption and the result follows by using Corollary 1 of \cite{Kounias1969} .\\
~\\
In the ARCH-type case when $f_\theta=0$ and
$A_0(h_\theta,\{\theta\})$ holds with $\tilde\beta^{(0)}(\theta)<1$,
we follow the same reasoning than previously remarking that
\eqref{ineqdep} has the simplified form:
$$
|q_s(\theta)-q_{s,j}(\theta)|\le \dfrac{1}{\underline{h}^2} X_s ^2
|h^{s}_\theta-h^{s,j}_\theta | + \dfrac{1}{\underline{h}} |
X_s^2-X_{s,j}^2 | + \dfrac{1}{\underline{h}}
|h^{s}_\theta-h^{s,j}_\theta |.
 $$
Then
$$
  (\E\|q_s(\theta)-q_{s,j}(\theta) \|_{{\Theta(r)}}^{r/4})^2\le C\E\big[\big( | X_s^2-X_{s,j}^2 |+\|h^{s}_\theta-h^{s,j}_\theta\|_{\Theta(r)}\big)^{r/2}\big].
$$
As
 $\|\|h^s_\theta - h^{s,j}_\theta \|_{\Theta(r)}\|_{r/2}\le C \sum\limits_{i\geq 1}\alpha_i^{(0)}(h_\theta,{\Theta(r)})\|X^2_{s-i}-X^2_{s-i,j}\|_{r/2}$ we derive from Lemma \ref{lem2},
 \begin{multline*}
 \E\ \|q_s(\theta)-q_{s,j}(\theta)\|_{{\Theta(r)}}^{r/4} \le C\Big[  \Big( \inf_{1\le p\le k/2}\big\{\tilde\beta^{(0)}(\theta_j^\ast)^{k/(2p)}
  +\sum_{i\ge p}\tilde\beta_i^{(0)}(\theta_j^\ast)\big\} \Big)^{r/4} +\Big(\sum\limits_{i\ge k/2}\alpha_i^{(0)}(h_\theta,{\Theta(r)})\Big)^{r/4}
\Big].
\end{multline*}
We easily conclude  to the result by choosing $p=u_k$ as above.\\
~\\
2-) We detail the proof for one order derivation in the general case where $A_0(f_\theta,\{\theta\})$ and $A_0(M_\theta,\{\theta\})$ hold with $\beta^{(0)}(\theta)<1$. The proofs of the other cases follow the same reasoning.\\
Let $j\in \{1,\cdots,K^*  \}$ and $i=1,\cdots,d$, we have:
   $$ \dfrac{v_{n^*_j}}{n^*_j}\Big\| \dfrac{\partial L_{n}(T^*_j, \theta)}{\partial \theta_i} - \dfrac{\partial L_{n,j}(T^*_j,\theta)}{\partial \theta_i}
 \Big\|_{{\Theta(r)}} \leq \frac{v_{n^*_j}}{n_j^*}\sum_{k=1}^{n_j^*}\Big\| \dfrac{\partial  q_{t^*_{j-1}+k}(\theta)}{\partial \theta_i}  -
    \dfrac{\partial  q_{t^*_{j-1}+k,j}(\theta)}{\partial \theta_i} \Big\|_{{\Theta(r)}}.$$
    By Corollary 1 of Kounias (1969), when $r\le 4$ with no loss of generality, it suffices to show
$$
  \sum_{k\ge 1} \big(\frac{v_k }{k}\big)^{r/4}  \E\Big( \Big\| \dfrac{\partial  q_{t^*_{j-1}+k}(\theta)}{\partial \theta_i}  -
    \dfrac{\partial  q_{t^*_{j-1}+k,j}(\theta)}{\partial \theta_i} \Big\|^{r/4}_{{\Theta(r)}} \Big)<\infty.
$$
For any $s\geq t^*_{j-1}$ denote $k=s- t^*_{j-1}$. For any $ \theta
\in {\Theta(r)}$, we have:
\begin{multline*}\dfrac{\partial q_{s}(\theta)}{\partial \theta_i} = -2\dfrac{(X_s-f^s_{\theta})}{h^s_{\theta}} \dfrac{\partial
  f^s_{\theta}}{\partial \theta_i}- \dfrac{(X_s-f^s_{\theta})^2}{(h^s_{\theta})^2} \dfrac{\partial
  h^s_{\theta}}{\partial \theta_i}+ \dfrac{1}{h^s_{\theta}} \dfrac{\partial h^s_{\theta}}{\partial \theta_i} \\
  \dfrac{\partial  q_{s,j}(\theta)}{\partial \theta_i} = -2\dfrac{(X_{s,j}-f^{s,j}_{\theta})}{h^{s,j}_{\theta}}
  \dfrac{\partial f^{s,j}_{\theta}}{\partial \theta_i}- \dfrac{(X_{s,j}-f^{s,j}_{\theta})^2}{(h^{s,j}_{\theta})^2}
  \dfrac{\partial h^{s,j}_{\theta}}{\partial \theta_i}+ \dfrac{1}{h^{s,j}_{\theta}} \dfrac{\partial
  h^{s,j}_{\theta}}{\partial \theta_i}.
\end{multline*}
Thus, using
  $|a_1 b_1 c_1 - a_2 b_2 c_2|\leq|a_1 - a_2||b_2||c_2| +  |b_1 - b_2||a_1||c_2| + |c_1 - c_2||a_1||b_1| $,
\begin{eqnarray}
 \nonumber  \Big \| \dfrac{\partial  q_{s}(\theta)}{\partial \theta_i} - \dfrac{\partial  q_{s,j}(\theta)}{\partial \theta_i} \Big \|_{{\Theta(r)}}
 \leq&& \hspace{-6mm}   2 \Big( \dfrac{1}{\underline{h}^2} \|h^s_{\theta}- h^{s,j}_{\theta}\|_{{\Theta(r)}} \|X_{s,j}-f^{s,j}_{\theta}\|_{{\Theta(r)}} \Big \|
  \dfrac{\partial f^{s,j}_{\theta}}{\partial \theta_i}\Big \|_{{\Theta(r)}}\\
 \nonumber  && \hspace{-3.5cm} +  \dfrac{1}{\underline{h}}(|X_s-X_{s,j}|+ \|f^s_{\theta}- f^{s,j}_{\theta}\|_{{\Theta(r)}})
   \Big \|\dfrac{\partial f^{s,j}_{\theta}}{\partial \theta_i}\Big \|_{{\Theta(r)}} + \dfrac{1}{\underline{h}} \Big \|\dfrac{\partial
    f^{s}_{\theta}}{\partial \theta_i} - \dfrac{\partial f^{s,j}_{\theta}}{\partial \theta_i}\Big\|_{{\Theta(r)}}
    \|X_{s}-f^{s}_{\theta} \|_{{\Theta(r)}} \Big)\\
  \nonumber  && \hspace{-3.5cm}  +  \dfrac{2}{\underline{h}^3} \|h^s_{\theta}- h^{s,j}_{\theta}\|_{{\Theta(r)}} \|X_{s,j}-f^{s,j}_{\theta} \|^2_{{\Theta(r)}}
        \Big \|\dfrac{\partial h^{s,j}_{\theta}}{\partial \theta_i}\Big \|_{{\Theta(r)}} \\
   \nonumber &&\hspace{-3.5cm}    + \dfrac{1}{\underline{h}} (|X_s-X_{s,j}|+ \|f^s_{\theta}- f^{s,j}_{\theta}\|_{{\Theta(r)}})(|X_s+X_{s,j}|+ \|f^s_{\theta}+
     f^{s,j}_{\theta}\|_{{\Theta(r)}}) \Big \|\dfrac{\partial f^{s,j}_{\theta}}{\partial \theta_i}\Big \|_{{\Theta(r)}}  \\
  \nonumber  &&\hspace{-3.5cm}  + \dfrac{1}{\underline{h}^2} \Big \|  \dfrac{\partial h^{s}_{\theta}}{\partial \theta_i}   -
        \dfrac{\partial h^{s,j}_{\theta}}{\partial \theta_i}\Big \|_{{\Theta(r)}} \|X_{s}-f^{s}_{\theta} \|^2_{{\Theta(r)}} +
     \dfrac{1}{\underline{h}^2} \|h^s_{\theta}- h^{s,j}_{\theta}\|_{{\Theta(r)}} \Big \| \dfrac{\partial h^{s,j}_{\theta}}{\partial \theta_i}\Big
     \|_{{\Theta(r)}}+  \dfrac{1}{\underline{h}} \Big \|  \dfrac{\partial h^{s}_{\theta}}{\partial \theta_i}   - \dfrac{\partial h^{s,j}_{\theta}}{\partial \theta_i}\Big \|_{{\Theta(r)}}
\end{eqnarray}
 So for all $s\ge t^*_{j-1}$ it holds:
\begin{multline}
 \big \| \dfrac{\partial  q_{s}(\theta)}{\partial \theta_i} - \dfrac{\partial  q_{s,j}(\theta)}{\partial \theta_i} \big \|_{{\Theta(r)}}\\
  \leq C \Big(  1+|X_s|^2+ |X_{s,j}|^2 +  \|f^s_{\theta}\|^2_{{\Theta(r)}} +  \|f^{s,j}_{\theta}\|^2_{{\Theta(r)}} +
   \big \| \dfrac{\partial f^{s}_{\theta}}{\partial \theta_i}\big \|^2_{{\Theta(r)}}+  \big \| \dfrac{\partial f^{s,j}_{\theta}}{\partial \theta_i}\big \|^2_{{\Theta(r)}}
   + \big \| \dfrac{\partial h^{s}_{\theta}}{\partial \theta_i}\big \|^2_{{\Theta(r)}}+  \big \| \dfrac{\partial h^{s,j}_{\theta}}{\partial \theta_i}\big \|^2_{{\Theta(r)}}
   \Big)\\
\quad \times \Big( |X_s-X_{s,j}| +  \|f^s_{\theta}-
f^{s,j}_{\theta}\|_{{\Theta(r)}} + \|h^s_{\theta}-
h^{s,j}_{\theta}\|_{{\Theta(r)}}  +
   \big \|  \dfrac{\partial f^{s}_{\theta}}{\partial \theta_i}   - \dfrac{\partial f^{s,j}_{\theta}}{\partial \theta_i}\big \|_{{\Theta(r)}}
  \nonumber + \Big \|  \dfrac{\partial h^{s}_{\theta}}{\partial \theta_i}   - \dfrac{\partial h^{s,j}_{\theta}}{\partial \theta_i}\Big \|_{{\Theta(r)}}\Big)
 \end{multline}
Since the processes admits finite moments of order $r$,  by
Cauchy-Schwartz Inequality:
\begin{multline}
 \Big( \E\big \| \dfrac{\partial  q_{s}(\theta)}{\partial \theta_i} - \dfrac{\partial  q_{s,j}(\theta)}{\partial \theta_i} \big
  \|^{r/4}_{{\Theta(r)}} \Big)^2 \leq C \Big( \E|X_s-X_{s,j}|^{r/2} + \E(\|f^s_{\theta}- f^{s,j}_{\theta}\|^{r/2}_{{\Theta(r)}})
   + \E(\|h^s_{\theta}- h^{s,j}_{\theta}\|^{r/2}_{{\Theta(r)}}) \\
   + \E\big \|\dfrac{\partial f^{s}_{\theta}}{\partial \theta_i} - \dfrac{\partial f^{s,j}_{\theta}}{\partial
  \theta_i}\big\|^{r/2}_{{\Theta(r)}} +  \E\big \|\dfrac{\partial h^{s}_{\theta}}{\partial \theta_i} - \dfrac{\partial h^{s,j}_{\theta}}{\partial
 \nonumber \theta_i}\big\|^{r/2}_{{\Theta(r)}}\Big)
 \end{multline}
As $(A_0(M_\theta,\Theta(r)))$ and $(A_1(M_\theta,\Theta(r)))$ hold
necessarily in this case, with the arguments of the proof of 1-),
for all $s\ge t^*_{j-1}$,
\begin{multline}
  \E\big \| \dfrac{\partial  q_{s}(\theta)}{\partial \theta_i} - \dfrac{\partial  q_{s,j}(\theta)}{\partial \theta_i} \big
  \|^{r/4}_{{\Theta(r)}}  \leq  C\Big[  \Big( \inf_{1\le p\le k/2}\big\{\beta^{(0)}(\theta_j^\ast)^{k/(2p)}
  +\sum_{i\ge p}\beta_i^{(0)}(\theta_j^\ast)\big\} \Big)^{r/4} +\Big(\sum\limits_{i\ge k/2}\alpha_i^{(0)}(f_\theta,{\Theta(r)})\Big)^{r/4}\\
   +\Big(\sum\limits_{i\ge k/2}\alpha_i^{(0)}(M_\theta,{\Theta(r)})\Big)^{r/4} +  \Big(\sum\limits_{i\ge k/2}\alpha_i^{(1)}(f_\theta,{\Theta(r)})\Big)^{r/4}
  \nonumber  + \Big(\sum\limits_{i\ge k/2}\alpha_i^{(1)}(M_\theta,{\Theta(r)})\Big)^{r/4}    \Big]
 \end{multline}
 Choosing  $p=u_k = kc^*/\log(k)$, we show (as in proof of 1-) ) that:
$$   \sum_{k\ge 1} \big(\frac{v_k }{k}\big)^{r/4}  \E\Big( \Big\| \dfrac{\partial  q_{t^*_{j-1}+k}(\theta)}{\partial \theta_i}  -
    \dfrac{\partial  q_{t^*_{j-1}+k,j}(\theta)}{\partial \theta_i} \Big\|^{r/4}_{{\Theta(r)}} \Big)<\infty.$$
\end{dem}
\vspace{-1cm}

\subsection{Consistency when the breaks are known}
When the breaks are known, we can chose $v_n =1$ for all $n$ in the
penalization of \eqref{contrast} as the penalization term does not
matter at all.
\begin{Prop}\label{prop1bis}
For all $j=1,\ldots,K^*$, under the assumptions of Lemma
\ref{lem_conv_L} 1-) with $v_n =1$ for all $n$, if the assumption
Id(${\Theta(r)}$) holds then
$$
\widehat{\theta}_n(T_j^*) \overset{\texttt{a.s.}}{\underset{n\to
\infty}\longrightarrow}\theta^*_j.
$$
\end{Prop}
\begin{dem} Let us first give the following useful corollary of Lemma \ref{lem_conv_L}
\begin{cor}\label{cor1}
\begin{itemize}
\item [i-)]  under the assumptions of Lemma \ref{lem_conv_L} 1-)  we have:\\
 $$\Big\|\dfrac{1}{n_j^*}\widehat{L}_n\big(T^*_j,\theta\big)-\mathcal{L}_{j}(\theta)\Big\|_{\Theta(r)}
 \overset{\texttt{a.s.}}{\underset{n\to \infty}\longrightarrow} 0~~\mbox{with}~~   \mathcal{L}_j(\theta) = -\dfrac{1}{2}\E\left(q_{0,j}(\theta)\right).$$
\item [ii-)] Under assumptions of Lemma \ref{lem_conv_L} 2-)   we have:
$$\Big\| \dfrac{1}{n^*_j} \dfrac{\partial^i \widehat{L}_n(T^*_j, \theta)}{\partial \theta^i} -
 \dfrac{\partial^i \mathcal{L}_j(\theta)}{\partial \theta^i} \Big\|_{{\Theta(r)}} \overset{\texttt{a.s.}}{\underset{n\to \infty}\longrightarrow} 0
~~\text{with} ~~  \dfrac{\partial^i \mathcal{L}_j(\theta)}{\partial
\theta^i}= -\dfrac{1}{2}\E\left(\dfrac{\partial^i
q_{0,j}(\theta)}{\partial
 \theta^i}\right).$$
\end{itemize}
\end{cor}
We conclude the proof of Proposition \ref{prop1bis} using
 $\mathcal{L}_{j}(\theta)= -\dfrac{1}{2}\E\left(q_{0,j}(\theta)\right)$ has a unique maximum in $\theta^*_j$ (see  \cite{Jeantheau1998}). From the almost sure convergence of the quasi-likelihood in i-) of Corollary \ref{cor1}, it comes:
  $$
   \widehat{\theta}_n(T_j^*)= \underset{\theta\in {\Theta(r)}} {\mbox{Argmax}} \left( \dfrac{1}{n_j^*}\widehat{L}_n \big(T^*_j,\theta\big) \right)
  \overset{\texttt{a.s.}}{\underset{n\to \infty}\longrightarrow}\theta^*_j.
  $$
\end{dem}
{\bf Proof of Corollary \ref{cor1}}
     Note that the proof of Lemma \ref{lem_conv_L} can be repeated by replacing $L_n$
     by the quasi-likelihood $\widehat{L}_n$. Thus, we obtain for $i=0,1,2$,
\begin{equation}\label{rem3}
 \frac{v_{n^\ast_j}}{n_j^*}\left\|\frac{\partial^i \widehat{L}_n\big(T^*_j,\theta\big)}{\partial
 \theta^i}- \frac{\partial^i L_{n,j}\big(T^*_j,\theta\big)}{\partial
 \theta^i}\right\|_{\Theta(r)}\underset{n\rightarrow\infty}{\longrightarrow}0.
\end{equation}
\begin{itemize}
 \item [i-)] Let $j\in {1,\cdots, K^*}$. From  \cite{Bardet2009}, we have:
$$
\Big\|\dfrac{1}{n_j^*}L_{n,j}\big(T^*_j,\theta\big)-\mathcal{L}_{j}(\theta)\Big\|_{\Theta(r)}
\overset{\texttt{a.s.}}{\underset{n\to \infty}\longrightarrow} 0.$$
Using \eqref{rem3}, the convergence to the limit likelihood follows.
\item [ii-)] From Lemma 4  and Theorem 1  of \cite{Bardet2009}, $\Big\| \dfrac{1}{n^*_j}
 \dfrac{\partial^i L_{n,j}(T^*_j, \theta)}{\partial \theta^i} - \dfrac{\partial^i \mathcal{L}_j(\theta)}{\partial \theta^i}
 \Big\|_{{\Theta(r)}} \overset{\texttt{a.s.}}{\underset{n\to \infty}\longrightarrow}
  0$ for $i = 1, \ 2$ and we conclude from \eqref{rem3}.\hfill$\blacksquare$
\end{itemize}

\subsection{Proof of Theorem \ref{theo4}} This proof is
 divided into two parts. In \textbf{part (1)} $K^*$ is assumed to be known and we show $(\widehat{\underline{\tau}}_n, \widehat{\underline{\theta}}_n) \limiteproban (\underline{\tau}^*,\underline{\theta}^*)$. In \textbf{part (2)}, $K^*$ is unknown and we show $ \widehat{K}_n
 \limiteproban K^*$ which ends the proof of Theorem  \ref{theo4}.\\

\textbf{Part (1)}.  Assume that $K^* $ is known and denote for any $
\underline{t}\in \mathcal{F}_{K^*}$:
\begin{eqnarray*}
 \widehat{I}_n(\underline{t}):=\widehat J_n(K^*,\underline{t},\underline{\widehat{\theta}}_n(\underline{t}))=
 -2\sum_{k=1}^{K^*} \sum_{j=1}^{K^*} \widehat L_n\left(T_k\cap T^*_j,\widehat{\theta}_n(T_k)\right)
\end{eqnarray*}
It comes that $\widehat{\underline{t}}_n= \underset{\underline{t}\in
\mathcal{F}_{K^*}}{\mbox{Argmin}}
     \left(\widehat{I}_n(\underline{t})\right)$. We show that   $ \widehat{\underline{\tau}}_n \overset{\texttt{P}}{\underset{n\to \infty}\longrightarrow} \underline{\tau}^*$  as it implies $ \widehat{\theta}_n(\widehat{T}_{n,j}) -  \widehat{\theta}_n(T_{j}^*)\overset{\texttt{P}}{\underset{n\to \infty}\longrightarrow} 0$ and from  Proposition \ref{prop1bis}
 $  \widehat{\theta}_n(\widehat{T}_{n,j})
 \overset{\texttt{P}}{\underset{n\to \infty}\longrightarrow} \theta_j^* $ for all $j=1, \cdots, K^*$. Without loss of generality, assume that $K^*=2$ and let $(u_n)$ be a sequence of positive integers satisfying
   $u_n\to\infty $, $ {u_n}/{n} \to0 $ and for some  $0 < \eta < 1$
\begin{eqnarray*}V_{\eta,u_n}&=& \{~t \in \Z / ~~ |t-t^*|>\eta n ~;~ u_n\leq t \leq n-u_n  ~\},\\
W_{\eta,u_n}&=& \{~t \in \Z / ~~ |t-t^*|>\eta n ~;~ 0<t<u_n ~~
\text{or} ~~ n-u_n<t\leq n  ~\}.
 \end{eqnarray*}
Asymptotically, we have $\mathbb{P}( \|\widehat{\underline{\tau}}_n
- \underline{\tau}^* \|_m >\eta )\simeq \mathbb{P}(|\widehat{t}_n -
t^*
 | >\eta n )$. But
\begin{eqnarray*}
\mathbb{P}(|\widehat{t}_n - t^*
 | >\eta n )  &\leq& \mathbb{P}\Big( \widehat{t}_n \in V_{\eta,u_n} \Big) + \mathbb{P}\Big( \widehat{t}_n \in W_{\eta,u_n} \Big) \\
    \nonumber  &\le& \mathbb{P}\Big( \underset{t\in V_{\eta,u_n}}{\mbox{min}}(\widehat{I}_n(t)- \widehat{I}_n(t^*)) \leq 0  \Big) +
                  \mathbb{P}\Big( \underset{t\in W_{\eta,u_n}}{\mbox{min}}(\widehat{I}_n(t)- \widehat{I}_n(t^*)) \leq 0  \Big)
\end{eqnarray*}
we show with similar arguments that these two probabilities tend to 0. We only detail below  the proof of $ \mathbb{P}\Big( \underset{t\in V_{\eta,u_n}}{\mbox{min}}(\widehat{I}_n(t)- \widehat{I}_n(t^*)) \leq 0  \Big)\to 0$ for shortness.\\
~\\
Let $t\in V_{\eta,u_n}$ satisfying $t^*\leq t$ (with no loss of
generality), then $T_1\cap T_1^*=T_1^*, ~~ T_2\cap T_1^*=\emptyset$
and $T_2\cap
T_2^*=T_2$. We decompose:%
\begin{multline}
\widehat{I}_n(t)-\widehat{I}_n(t^*)=
2\Big(\widehat{L}_n(T^*_1,\widehat{\theta}_n(T^*_1))-\widehat{L}_n(T^*_1,\widehat{\theta}_n(T_1))+
   \widehat{L}_n(T_1\cap T^*_2,\widehat{\theta}_n(T_2^*))\\
 \label{diff_I} -\widehat{L}_n(T_1\cap T^*_2,\widehat{\theta}_n(T_1))+\widehat{L}_n(T_2,\widehat{\theta}_n(T_2^*))-\widehat{L}_n(T_2,\widehat{\theta}_n(T_2))\Big).
\end{multline}
As  $\# T_1^*=t^*, ~\#(T_1\cap T_2^*)=t-t^*, ~\#T_2 =n- t\geq u_n$,
each  term tends  to $\infty$ with $n$. Using Proposition
\ref{prop1bis} and Corollary \ref{cor1}, we get the following
convergence, uniformly on  $V_{\eta,u_n}$,
\begin{multline*} \widehat{\theta}_n(T_1^*) \limitepsn \theta^*_1,~~
 \widehat{\theta}_n(T_2^*) \limitepsn\theta^*_2,~~\widehat{\theta}_n(T_2) \limitepsn\theta^*_2~~\mbox{ and }~~
 \Big\|\dfrac{\widehat{L}_n\big(T_1^*,\theta\big)}{n}-\tau_1^*\mathcal{L}_{1}(\theta)\Big\|_{\Theta(r)}
  \limitepsn0,\\
\Big\|\dfrac{\widehat{L}_n\big(T_1\cap
T_2^*,\theta\big)}{t-t^*}-\mathcal{L}_{2}(\theta)\Big\|_{\Theta(r)}
 \limitepsn 0,~~
 \Big\|\dfrac{\widehat{L}_n\big(T_2,\theta\big)}{n-t}-\mathcal{L}_{2}(\theta)\Big\|_{\Theta(r)}
 \limitepsn0.\end{multline*}
For any $\varepsilon>0$, there exists an integer $N_0$ such that for
any  $n>N_0$,
$$\Big\|\dfrac{\widehat{L}_n\big(T_1^*,\theta\big)}{n}-\tau_1^*\mathcal{L}_{1}(\theta)\Big\|_{\Theta(r)}<\dfrac{\varepsilon}{6}
; ~~
 \Big\|\dfrac{\widehat{L}_n\big(T_1\cap
T_2^*,\theta\big)}{t-t^*}-\mathcal{L}_{2}(\theta)\Big\|_{\Theta(r)}<\dfrac{\varepsilon}{6};
~~
\Big|\dfrac{\widehat{L}_n\big(T_1^*,\widehat{\theta}_n(T_1^*)\big)}{n}-\tau_1^*\mathcal{L}_{1}(\theta_1^*)\Big|<\dfrac{\varepsilon}{6}$$
 $$  \Big|\dfrac{\widehat{L}_n\big(T_1\cap
   T_2^*,\widehat{\theta}_n(T_2^*)\big)}{t-t^*}-\mathcal{L}_{2}(\theta_2^*)\Big|<\dfrac{\varepsilon}{6}; ~~~
\dfrac{n-t}{n}\Big|\dfrac{\widehat{L}_n(T_2,\widehat{\theta}_n(T_2^*))-\widehat{L}_n(T_2,\widehat{\theta}_n(T_2))}{n-t}\Big|<\dfrac{\varepsilon}{6}$$
Thus, for $n>N_0$,
\begin{eqnarray}
\nonumber  \tau_1^*\mathcal{L}_1(\theta_1^*)-
\tau_1^*\mathcal{L}_1(\widehat{\theta}_n(T_1))&=&
\tau_1^*\mathcal{L}_1(\theta_1^*)-
\dfrac{\widehat{L}_n\big(T_1^*,\widehat\theta_n(T_1^*)\big)}{n}
     + \dfrac{\widehat{L}_n\big(T_1^*,\widehat{\theta}_n(T_1^*)\big)}{n} -
        \dfrac{\widehat{L}_n\big(T_1^*,\widehat{\theta}_n(T_1)\big)}{n}\\
\nonumber  && \hspace{5cm}  + \dfrac{\widehat{L}_n\big(T_1^*,\widehat{\theta}_n(T_1)\big)}{n} -  \tau_1^*\mathcal{L}_1(\widehat{\theta}_n(T_1))\\
  \nonumber  &\leq& \dfrac{\varepsilon}{6} + \dfrac{\widehat{L}_n\big(T_1^*,\widehat{\theta}_n(T_1^*)\big)}{n} - \dfrac{\widehat{L}_n\big(T_1^*,\widehat{\theta}_n(T_1)\big)}{n}+
  \dfrac{\varepsilon}{6}.
\end{eqnarray}
Then,
\begin{eqnarray} \label{diff_I1}
    \dfrac{\widehat{L}_n\big(T_1^*,\widehat{\theta}_n(T_1^*)\big)}{n} - \dfrac{\widehat{L}_n\big(T_1^*,\widehat{\theta}_n(T_1)\big)}{n}&>&
            \tau_1^*\Big(\mathcal{L}_1(\theta_1^*)- \mathcal{L}_1(\widehat{\theta}_n(T_1))\Big)- \dfrac{\varepsilon}{3}.
\end{eqnarray}
Similarly, for $n>N_0$:
\begin{eqnarray} \label{diff_I2}
   ~~~~~~~~ \dfrac{\widehat{L}_n\big(T_1\cap T_2^*,\widehat{\theta}_n(T_2^*)\big)}{n} - \dfrac{\widehat{L}_n\big(T_1\cap T_2^*,\widehat{\theta}_n(T_1)\big)}{n}>
            \eta \Big( \mathcal{L}_{2}(\theta_2^*)-  \mathcal{L}_{2}(\widehat{\theta}_n(T_1))\Big)- \dfrac{\varepsilon}{3}.
\end{eqnarray}
Finally,  for $n>N_0$,
\begin{eqnarray} \label{diff_I3}
    \dfrac{\widehat{L}_n(T_2,\widehat{\theta}_n(T_2^*))-\widehat{L}_n(T_2,\widehat{\theta}_n(T_2))}{n}&>&
            - \dfrac{\varepsilon}{6},
\end{eqnarray}
and from  \eqref{diff_I} and inequalities \eqref{diff_I1},
\eqref{diff_I2} and \eqref{diff_I3} we obtain uniformly in $t$:
$$ \frac{ \widehat{I}_n(t)-\widehat{I}_n(t^*)}{n}> \tau_1^*\Big(\mathcal{L}_1(\theta_1^*)- \mathcal{L}_1(\widehat{\theta}_n(T_1))\Big)+
        \eta \Big( \mathcal{L}_{2}(\theta_2^*)-  \mathcal{L}_{2}(\widehat{\theta}_n(T_1))\Big)- \dfrac{5}{6}\varepsilon, \quad n>N_0.$$
Since $\theta_1^* \neq \theta_2^*$, let $\mathcal{V}_1$,
$\mathcal{V}_2$ be  two open neighborhoods and disjoint of
$\theta_1^*$ and  $\theta_2^*$ respectively,
  $$
  \delta_i:=\underset{\theta\in\mathcal{V}^{c}_i} {\mbox{Inf}}\Big( \mathcal{L}_i(\theta_i^*)-\mathcal{L}_i(\theta) \Big)>0~~\mbox{ for  }~~i=1,2,
$$
since the  function $\theta \mapsto \mathcal{L}_j(\theta)$ has a
strict maximum in $\theta_j^*$  (see \cite{Jeantheau1998}). With
$\varepsilon=\min(\tau_1^*\delta_1,~\eta \delta_2)$, we get
\begin{itemize}
 \item if $\widehat{\theta}_n(T_1) \in \mathcal{V}_1$ i.e. $\widehat{\theta}_n(T_1) \in \mathcal{V}_2^c$, then
  $\dfrac{\widehat{I}_n(t)-\widehat{I}_n(t^*)}{n}> \eta \delta_2- \dfrac{5}{6}\varepsilon \geq \dfrac{\varepsilon}{6}$;
 \item If $\widehat{\theta}_n(T_1) \notin \mathcal{V}_1$ i.e. $\widehat{\theta}_n(T_1) \in \mathcal{V}_1^c$, then
  $\dfrac{\widehat{I}_n(t)-\widehat{I}_n(t^*)}{n}>\tau_1^* \delta_1- \dfrac{5}{6}\varepsilon \geq \dfrac{\varepsilon}{6}$.
\end{itemize}
In any case we prove that $\widehat{I}_n(t)-\widehat{I}_n(t^*) >
\dfrac{\varepsilon}{6}n $ for  $n>N_0$ and all $t \in V_{\eta,
u_n}$. It implies that $ \mathbb{P}\Big( \underset{\underline{t}\in
V_{\eta,u_n}}{\mbox{min}}(\widehat{I}_n(t)- \widehat{I}_n(t^*)) \leq
0 \Big){\underset{n\to \infty}\longrightarrow}0$ and we show
similarly  $ \mathbb{P}\Big( \underset{\underline{t}\in
W_{\eta,u_n}}{\mbox{min}}(\widehat{I}_n(t)- \widehat{I}_n(t^*)) \leq
  0 \Big){\underset{n\to \infty}\longrightarrow}0. $ It follows directly that
 $\mathbb{P} (\|\widehat{\underline{\tau}}_n - \underline{\tau}^* \|_m >\eta  ) {\underset{n\to \infty}\longrightarrow} 0 ~
\text{for all }~ \eta > 0. $\\

\textbf{Part(2)}. Now $K^*$ is unknown. For  $K\geq 2$,
$x=(x_1,\cdots,x_{K-1}) \in \R^{K-1}$,  $y=(y_1,\cdots,y_{K^*-1})
  \in \R^{K^*-1}$, denote
  $$\| x-y \|_{\infty} = \underset{1\leq j\leq K^*-1} {\mbox{max}} ~~ \underset{1\leq k\leq K-1} {\mbox{min}} | x_k-y_j | .$$ The following Lemma follows directly from
\textbf{Part(1)} and the definition of $\|\cdot\|_\infty$:
\begin{lem} \label{lem4}
 Let $K\geq 1$, $(\widehat{\underline{t}}_n, \widehat{\underline{\theta}}_n)$
  obtained by the minimization of $\widehat{J}_n(\underline{t}, \underline{\theta})$ on $\mathcal{F}_{K} \times {\Theta(r)}^K$
  and $\widehat{\underline{\tau}}_n = {\widehat{\underline{t}}_n}/{n}$. Under assumptions of Theorem \ref{theo4},
   $\| \widehat{\underline{\tau}}_n - \underline{\tau}^* \|_{\infty} \overset{\texttt{P}}{\underset{n\to
  +\infty}\longrightarrow} 0 $  if  $ K \geq  K^*$.
\end{lem}
Now we use the following Lemma \ref{lem5} which is proved below (see
also \cite{LavielleMoulines2000}):
\begin{lem}\label{lem5} Under the assumptions of Lemma \ref{lem_conv_L} i-), for any $K \geq 2$, there exists $C_K > 0$ such as:
 $$ \forall  (\underline{t}, \underline{\theta}) \in \mathcal{F}_{K} \times {\Theta(r)}^K,  ~~ u_n (\underline{t}, \underline{\theta})=
  2 \sum_{j=1}^{K^*} \sum_{k=1}^{K} \dfrac{n_{kj}}{n} ( \mathcal{L}_j(\theta^*_j)-
  \mathcal{L}_j(\theta_k)) \geq \dfrac{C_K}{n} \| \underline{t}- \underline{t}^* \|_{\infty}. $$
\end{lem}
Continue with the proof of \textbf{Part(2)} shared in two parts,
{\it i.e.} we show that $ \text{P}( \widehat{K}_n = K
)\overset{\texttt{}}{\underset{n\to + \infty}\longrightarrow} 0 $
for $K<K^*$ and $K^*<K\le K_{\text{max}}$ separately. In any case,
we have \begin{eqnarray}
  \text{P}(\widehat{K}_n = K)
   \nonumber & \leq & \text{P}\Big( \underset{( \underline{t}, \underline{\theta}) \in \mathcal{F}_{K} \times {\Theta(r)}^K } {\mbox{inf}}
   ( \widetilde{J}_n ( K, \underline{t}, \underline{\theta} )) \leq \widetilde{J}_n (K^*,\underline{t}^*,\underline{\theta}^*)
   \Big)\\
   \label{k_inf1}  & \leq & \text{P}\Big( \underset{( \underline{t}, \underline{\theta}) \in \mathcal{F}_{K} \times {\Theta(r)^K} } {\mbox{inf}}
   ( \widehat{J}_n (K, \underline{t}, \underline{\theta} ) - \widehat{J}_n (K^*,\underline{t}^*, \underline{\theta}^* ) ) \leq \frac{n}{v_n} (K^*-K).
   \Big).
   \end{eqnarray}

\begin{itemize}
\item [i-)] For $K<K^*$, we decompose $ \widehat{J}_n (K, \underline{t},
  \underline{\theta} ) - \widehat{J}_n (K^*,\underline{t}^*, \underline{\theta}^*)  = n(u_n (\underline{t}, \underline{\theta})
   + e_n (\underline{t}, \underline{\theta})) $  where $u_n$ is defined in Lemma \ref{lem5} and
$$
 e_n (\underline{t}, \underline{\theta}) =  2 \left[  \sum_{j=1}^{K^*} \dfrac{n^*_j}{n} \Big(
 \dfrac{\widehat{L}_n(T^*_j,  \theta^*_j)}{n^*_j} - \mathcal{L}_j(\theta^*_j) \Big) +  \sum_{k=1}^{K} \sum_{j=1}^{K^*} \dfrac{n_{kj}}{n}
 \Big(  \mathcal{L}_j(\theta_k) - \dfrac{\widehat{L}_n(T^*_j\cap T_k, \theta_k)}{n_{kj}} \Big) \right] .
$$
It comes from the relation \eqref{k_inf1} that:
\begin{eqnarray}
  \text{P}(\widehat{K}_n = K)  \leq
   \label{k_inf2} \text{P}\Big( \underset{( \underline{t}, \underline{\theta}) \in \mathcal{F}_{K} \times \Theta_K } {\mbox{inf}}
    (u_n (\underline{t}, \underline{\theta}) + e_n (\underline{t}, \underline{\theta})) \leq \dfrac{\beta_n}{n} (K^*-K)
   \Big).
   \end{eqnarray}
Corollary \ref{cor1} ensures that  $ e_n
(\underline{t},\underline{\theta})) \to 0$ a.s. and uniformly  on
$\mathcal{F}_{K} \times \Theta(r)^K$. By Lemma \ref{lem5}, there
exists $C_K >0$ such that
    $u_n (\underline{t}, \underline{\theta}) \geq {C_K} \| \underline{t}- \underline{t}^* \|_{\infty}/{n}$ for all
     $  (\underline{t}, \underline{\theta}) \in \mathcal{F}_{K} \times \Theta(r)^K $.
      But, since $K<K^*$, for any  $\underline{t} \in \mathcal{F}_{K}$, we have $ \| \underline{t}- \underline{t}^* \|_{\infty}/n=
      \| \underline{\tau}- \underline{\tau}^* \|_{\infty} \geq  \min_{1 \leq j \leq K^* } ( \tau^*_j - \tau^*_{j-1} )/2$ that is positive by assumption. Then $ u_n (\underline{t}, \underline{\theta}) >0$ for  all
     $  (\underline{t}, \underline{\theta}) \in \mathcal{F}_{K} \times \Theta(r)^K $ and since $ 1/v_n \overset{\texttt{}}{\underset{n\to \infty}\longrightarrow} 0
    $, we deduce from \eqref{k_inf2} that $\text{P}(\widehat{K}_n = K) \overset{\texttt{}}{\underset{n\to \infty}\longrightarrow} 0 .$
\item[ii-)] Now let $K^*<K\le K_{\text{max}}$. from \eqref{k_inf2} and the Markov Inequality we have:
\begin{eqnarray}
  \text{P}(\widehat{K}_n = K)
   \nonumber
   \nonumber & \leq & \text{P} \Big( \widehat{J}_n (K, \widehat{\underline{t}}_n, \widehat{\underline{\theta}}_n ) - \widehat{J}_n (K^*,\underline{t}^*, \underline{\theta}^* ) + \frac{n}{v_n}
   (K-K^*) \leq 0  \Big)\\
   \nonumber & \leq & \text{P} \Big( | \widehat{J}_n (K, \widehat{\underline{t}}_n, \widehat{\underline{\theta}}_n ) - \widehat{J}_n (K^*,\underline{t}^*, \underline{\theta}^* ) | \geq\frac{n}{v_n}
         \Big)\\
      \label{k_inf3} & \leq & \frac{v_n}{n}\mathbb{E}| \widehat{J}_n (K,\widehat{\underline{t}}_n, \widehat{\underline{\theta}}_n ) -
  \widehat{J}_n (K^*,\underline{t}^*, \underline{\theta}^* ) |.
   \end{eqnarray}
Denote $ \widehat{\underline{t}}_n = (\widehat{t}_{n,1}, \cdots,
\widehat{t}_{n,K}) $. By  Lemma \ref{lem4}, there exists some subset
$ \{ k_j ~, 1
 \leq j \leq K^*-1  \} $ of $ \{ 1, \cdots,  K-1  \}  $ such that for any $ j = 1, \cdots,  K^*-1 $, $  {\widehat{t}_{n,k_j}}/{n}
 \overset{\texttt{}}\to \tau^*_j $. Denoting $k_0 = 0$  and  $  k_{K^*} =K$, we have:
\begin{eqnarray}
  \widehat{J}_n (K,\widehat{\underline{t}}_n, \widehat{\underline{\theta}}_n ) - \widehat{J}_n (K^*,\underline{t}^*, \underline{\theta}^* )
   \nonumber & = & 2 \Big( \sum_{j=1}^{K^*} \widehat{L}_n(T^*_j, \theta^*_j)-  \sum_{k=1}^{K} \widehat{L}_n( \widehat{T}_{n,k}, \widehat{\theta}_{n,k}) \Big)\\
  \nonumber  & = & 2 \sum_{j=1}^{K^*} \Big[ \widehat{L}_n(T^*_j,\theta^*_j) - \sum_{k=k_{j-1}+1}^{k_j} \widehat{L}_n( \widehat{T}_{n,k}, \widehat{\theta}_{n,k}) \Big]
   \end{eqnarray}
and from \eqref{k_inf3} we deduce that:
\begin{eqnarray}
   \text{P}(\widehat{K}_n = K)
   \nonumber & \leq & \dfrac{2v_n}{n} \sum_{j=1}^{K^*} \mathbb{E} \Big| \widehat{L}_n(T^*_j,\theta^*_j) - \sum_{k=k_{j-1}+1}^{k_j} \widehat{L}_n( \widehat{T}_{n,k}, \widehat{\theta}_{n,k}) \Big |  \\
  \nonumber  & \leq &  C\sum_{j=1}^{K^*} \dfrac{v_{n^*_j} }{n^*_j} \mathbb{E} \Big| \widehat{L}_n(T^*_j,\theta^*_j) -  \sum_{k=k_{j-1}+1}^{k_j} \widehat{L}_n( \widehat{T}_{n,k}, \widehat{\theta}_{n,k})
  \Big|.
   \end{eqnarray}
  Since for any $ j = 1, \cdots,  K^*-1$, it comes from Lemma \ref{lem_conv_L} that
   $$\frac{v_{n^*_j}}{n^*_j}  \mathbb{E} \Big|  \widehat{L}_n(T^*_j,\theta^*_j) - \sum_{k=k_{j-1}+1}^{k_j} \widehat{L}_n( \widehat{T}_{n,k},
   \widehat{\theta}_{n,k}) \Big |
    \overset{\texttt{}}{\underset{n\to \infty}\longrightarrow} 0, $$
and therefore  $ \text{P}(\widehat{K}_n = K)
\overset{\texttt{}}{\underset{n\to \infty}\longrightarrow} 0
$.~~~~~$\blacksquare$
\end{itemize}
{\bf Proof of Lemma \ref{lem5} } Let $K \geq 1$ and consider the
real function $\upsilon$ define on $\Theta \times \Theta$ by:
$$ \upsilon (\theta, \theta ') = \left\{ \begin{array}{l}
 \underset{ 1 \leq j \leq K^* } {\mbox{min}}[ \text{max}( \mathcal{L}_j(\theta^*_j) - \mathcal{L}_j(\theta),  ~ \mathcal{L}_j(\theta^*_j) - \mathcal{L}_j(\theta ') )  ]  ~~ \mbox{if} ~~ \theta  \ne \theta '  \\
 0 ~~\mbox{if} ~~ \theta  = \theta '. \\
 \end{array} \right. $$
The function $\upsilon$ has positive values and $ \upsilon (\theta
,\theta ') =0 $ if and only if $ \theta =\theta ' $ since the
function $\theta \mapsto \mathcal{L}_j(\theta)$ has a strict maximum
in $\theta_j^*$  (see \cite{Jeantheau1998}). By Lemma 3.3
 of \cite{Lavielle2000}, there exists  $C_{\theta^*} > 0$ such that for any
 $  (\underline{t}, \underline{\theta}) \in \mathcal{F}_{K} \times \Theta_K $
 $$ \sum_{j=1}^{K^*} \sum_{k=1}^{K} \dfrac{n_{kj}}{n}  \upsilon(\theta_k,  \theta^*_j) \geq \dfrac{C_{\theta^*}}{n} \| \underline{t}- \underline{t}^*
    \|_{\infty} .$$
  Moreover, for any  $ j = 1, \cdots, K^*$  and  $\theta \in \Theta$,  $  \mathcal{L}_j(\theta^*_j) - \mathcal{L}_j(\theta) \geq
  \upsilon(\theta,  \theta^*_j)$ and denoting $C_K = 2 C_{\theta^*}$ the result follows immediately. $\blacksquare$

\subsection{ Proof of Theorem \ref{theo3}}
Assume  with no loss of generality that $K^*=2$. Denote $ (u_n)_n $
a sequence  satisfying $u_n
   {\underset{n\to \infty}\longrightarrow} \infty $, ${u_n}/{n} {\underset{n\to \infty}\longrightarrow} 0 $ ~
  and ~ $ \mathbb{P}( | \widehat{\underline{t}}_n - \underline{t}^* | > u_n ) {\underset{n\to \infty}\longrightarrow} 0 $ (for example
   $u_n=  n \sqrt{\text{max}(\mathbb{E}| \widehat{\tau}_n - \tau^* |,  n^{-1})}$). For $\delta > 0$, as we have
$$ \mathbb{P}( | \widehat{\underline{t}}_n - \underline{t}^* | > \delta )
    \leq
    \mathbb{P}(  \delta < | \widehat{\underline{t}}_n - \underline{t}^* | \leq  u_n ) + \mathbb{P}( | \widehat{\underline{t}}_n - \underline{t}^* |_m > u_n )
$$
it suffices to show that $ \mathop {\lim }\limits_{\delta  \to \,
\infty } ~ \mathop {\lim }\limits_{ n\to
         \infty } \mathbb{P}( \delta <| \widehat{t}_n - t^*| \leq u_n  ) = 0$.\\
Denote $V_{\delta,u_n}= \{~t \in \Z / ~~ \delta < |t-t^*| \leq u_n
~\}$. Then,
   $$  \mathbb{P}( \delta <| \widehat{t}_n - t^*| \leq u_n  ) \leq  \mathbb{P}\Big( \underset{t\in V_{\delta,u_n}}{\mbox{min}}(\widehat{I}_n(t)- \widehat{I}_n(t^*)) \leq 0 \Big).   $$
 Let $t \in V_{\delta,u_n} $ (for example  $t\geq t^*$). With the notation of the proof of Theorem \ref{theo4}, we have
   $\widehat{L}_n(T^*_1,\widehat{\theta}_n(T^*_1)) \geq \widehat{L}_n(T^*_1,\widehat{\theta}_n(T_1))$ and from \eqref{diff_I} we obtain:
   $$  \dfrac{\widehat{I}_n(t)-\widehat{I}_n(t^*)}{t-t^*} \geq
  \dfrac{2}{t-t^*} \Big( \widehat{L}_n(T_1\cap T^*_2,\widehat{\theta}_n(T_2^*))-\widehat{L}_n(T_1\cap T^*_2,\widehat{\theta}_n(T_1))
  +\widehat{L}_n(T_2,\widehat{\theta}_n(T_2^*))-\widehat{L}_n(T_2,\widehat{\theta}_n(T_2))\Big) .  $$
We conclude in two steps:
\begin{description}
 \item [i-)] We show that  $   \dfrac{1}{t-t^*} \Big( \widehat{L}_n(T_1\cap T^*_2,\widehat{\theta}_n(T_2^*))-\widehat{L}_n(T_1\cap
      T^*_2,\widehat{\theta}_n(T_1)) \Big)>0$ for $n$ large enough. Then $ \dfrac{\widehat{L}_n(T_1, \theta)}{n} = \dfrac{t^*}{n} \dfrac{\widehat{L}_n(T^*_1, \theta)}{t^*}
 + \dfrac{t-t^*}{n} \dfrac{\widehat{L}_n(T_1 \cap T^*_2,\theta)}{t-t^*}$ and since  $\dfrac{t-t^*}{n} \leq \dfrac{u_n}{n} ~~ {\underset{n\to \infty}\longrightarrow} 0 $ and
$$
\widehat{\theta}_n(T_1)= \underset{\theta\in {\Theta(r)}}
{\mbox{Argmax}} \left( \dfrac{1}{n}\widehat{L}_n
    \big(T_1,\theta\big) \right) \limitepsu \theta_1^*.
$$
It comes that
  $   \dfrac{1}{t-t^*} \Big( \widehat{L}_n(T_1\cap T^*_2,\widehat{\theta}_n(T_2^*))-\widehat{L}_n(T_1\cap
      T^*_2,\widehat{\theta}_n(T_1)) \Big)$ converges a.s. and
      uniformly on  $ V_{\delta,u_n} $  to  $  \mathcal{L}_{2}(\theta^*_2)-  \mathcal{L}_{2}(\theta^*_1)>0.  $

 \item [ii-)] We show that  $ \dfrac{1}{t-t^*} \Big(\widehat{L}_n(T_2,\widehat{\theta}_n(T_2^*))-\widehat{L}_n(T_2,\widehat{\theta}_n(T_2))\Big)   \limitepsu 0$. For large value of $n$, we remark that  $ \widehat{\theta}_n(T_2) \in  \mathop {\Theta(r)} \limits^\circ$  so that $  {\partial \widehat{L}_{n}(T_2,
 \widehat{\theta}_n(T_2))}/{\partial \theta}  = 0 $. The mean value theorem on ${\partial \widehat{L}_{n}}/{\partial
 \theta_i}$ for any $i=1,\dots,d$ gives the existence of  $ \widetilde{\theta}_{n,i} \in [ \widehat{\theta}_n(T_2),
 \widehat{\theta}_n(T_2^*) ] $ such that:
\begin{eqnarray}
   0 & =&  \dfrac{\partial \widehat{L}_{n}(T_2,\widehat{\theta}_n(T^*_2))}{\partial \theta_i} +
    \dfrac{\partial^2 \widehat{L}_n(T_2, \widetilde{\theta}_{n,i})}{\partial \theta \partial
   \label{accr_f}  \theta_i} (\widehat{\theta}_n(T_2)) -  \widehat{\theta}_n(T_2^*))
   \end{eqnarray}
 where for  $a,b \in \mathbb{R}^d ~,~ [a,  b ] = \{ (1-\lambda)a + \lambda b ~ ; ~ \lambda \in [0,1 ] \}. $  Using the equalities $ \widehat{L}_{n}(T^*_2, \theta) = \widehat{L}_{n}(T_1 \cap T^*_2, \theta) + \widehat{L}_{n}(T_2,
  \theta)$  and   $ {\partial \widehat{L}_{n}(T^*_2, \widehat{\theta}_n(T^*_2))}/{\partial \theta}
  = 0 $,  it comes from \eqref{accr_f}:
$$ \dfrac{\partial \widehat{L}_{n}(T_1 \cap T^*_2,\widehat{\theta}_n(T^*_2))}{\partial \theta_i}
     =   \dfrac{\partial^2 \widehat{L}_n(T_2, \widetilde{\theta}_{n,i})}{\partial \theta \partial \theta_i} (\widehat{\theta}_n(T_2)) -  \widehat{\theta}_n(T_2^*)),\qquad \forall i=1,\dots,d,
     $$  and it follows:
\begin{equation}
  \dfrac{1}{t-t^*} \dfrac{\partial \widehat{L}_{n}(T_1 \cap T^*_2,\widehat{\theta}_n(T^*_2))}{\partial \theta}  =
   \label{accr_f2} \frac{n-t}{t-t^*}A_n \cdot(\widehat{\theta}_n(T_2) -  \widehat{\theta}_n(T_2^*))
\end{equation}
 with $ A_n:= \Big( \dfrac{1}{n-t} \dfrac{\partial^2 \widehat{L}_n(T_2, \widetilde{\theta}_{n,i})}{\partial \theta \partial \theta_i} \Big)_{1\leq i \leq d}.$  Corollary \ref{cor1} ii-) gives that:
$$
\dfrac{1}{t-t^*} \dfrac{\partial \widehat{L}_{n}(T_1 \cap
T^*_2,\widehat{\theta}_n(T^*_2))}{\partial \theta}
   \limitepsu \dfrac{\partial \mathcal{L}_{2}(\theta^*_2)}{\partial \theta} = 0
$$
 and $ A_n \limitepsu - \dfrac{1}{2}\mathbb{E}\Big( \dfrac{\partial^2 q_{0,2}(\theta^*_2)}{\partial \theta^2} \Big).$ Under assumption
 (Var),  $\mathbb{E}\Big( \dfrac{\partial^2 q_{0,2}(\theta^*_2)}{\partial \theta^2} \Big) $ is a nonsingular
 matrix (see \cite{Bardet2009}). Then, we deduce from
\eqref{accr_f2} that
\begin{equation}\label{vitesse_theta}
  \frac{n-t}{t-t^*}(\widehat{\theta}_n(T_2) -  \widehat{\theta}_n(T_2^*)) \limitepsu 0 .
  \end{equation}
We conclude by the Taylor expansion on $ \widehat{L}_n$ that gives
\begin{multline*}
 \dfrac{1}{t-t^*}     | \widehat{L}_n(T_2,\widehat{\theta}_n(T_2))-\widehat{L}_n(T_2,\widehat{\theta}_n(T_2^*))|\\
      \leq  \dfrac{1}{2(t-t^*)} \| \widehat{\theta}_n(T_2)) -  \widehat{\theta}_n(T_2^*)  \|^2
       \underset{\theta\in {\Theta(r)}} \sup \Big \|
           \label{f_Taylor}  \dfrac{\partial^2 \widehat{L}_n(T_2, \theta)}{\partial \theta^2} \Big \|\to 0\quad \mbox{a.s.}\qquad\blacksquare
\end{multline*}
\end{description}

 \subsection{ Proof of Theorem \ref{TLC}} First, $\big( \widehat{\theta}_n(\widehat{T}_j) -
\theta^*_j \big) = \big(
\widehat{\theta}_n(\widehat{T}_j) - \widehat{\theta}_n(T^*_j)\big) +
 \big( \widehat{\theta}_n(T^*_j) - \theta^*_j \big) $
for any $j\in \{1,\cdots,K^* \}$.
   By  Theorem \ref{theo3} it comes $\widehat{t}_j-t^*_j =
   o_P(\log(n))$. Using relation \eqref{vitesse_theta}, we obtain: $ \widehat{\theta}_n(\widehat{T}_j) - \widehat{\theta}_n(T^*_j)=
   o_P(\dfrac{\log(n)}{n})$.
   Hence, $ \sqrt{n^*_j} \big( \widehat{\theta}_n(\widehat{T}_j) - \widehat{\theta}_n(T^*_j)\big) \overset{P}{\underset{n\to \infty} \longrightarrow} 0$ and
it suffices to show that   $ \sqrt{n^*_j}
\big(\widehat{\theta}_n(T^*_j) - \theta^*_j \big) \overset{\mathcal{
D}}{\underset{n\to \infty} \longrightarrow}
   \mathcal{N}_d \big(0, F(\theta^*_j)^{-1}G(\theta^*_j)F(\theta^*_j)^{-1} \big)$ to
   conclude.\\
For large value of $n$,  $ \widehat{\theta}_n(T^*_j) \in \mathop
{\Theta(r)}\limits^\circ$. By  the mean value theorem, there exists
$ (\widetilde{\theta}_{n,k})_{1\leq k\leq d} \in [
\widehat{\theta}_n(T^*_j), \theta^*_j ]$ such that
\begin{equation}\label{eq_AN_accr_f}
 \dfrac{\partial L_{n}(T^*_j,\widehat{\theta}_n(T^*_j))}{\partial \theta_k}  =  \dfrac{\partial L_{n}(T^*_j,\theta^*_j)}{\partial \theta_k}
  + \dfrac{\partial^2 L_n(T^*_j, \widetilde{\theta}_{n,k})}{\partial \theta \partial \theta_k} (\widehat{\theta}_n(T^*_j) -
  \theta^*_j).
 \end{equation}
Let $ F_n= -2\Big(  \dfrac{1}{n^*_j} \dfrac{\partial^2 L_n(T^*_j,
\widetilde{\theta}_{n,k})}{\partial \theta \partial \theta_k}
\Big)_{1\leq k \leq d}.$ By  Lemma \ref{lem_conv_L} and Corollary
\ref{cor1}, $F_n \overset{\text{a.s.}}{\underset{n\to \infty}
\longrightarrow} F(\theta^*_j)$ (where $F(\theta^*_j)$ is defined by
\eqref{eq_AN}). But, under (Var), $F(\theta^*_j)$ is a non singular
matrix (see \cite{Bardet2009}). Thus, for $n$ large enough, $ F_n $
is invertible and \eqref{eq_AN_accr_f} gives
  $$  \sqrt{n^*_j} \big(\widehat{\theta}_n(T^*_j) - \theta^*_j \big) = -2F_n^{-1} \Big[ \dfrac{1}{\sqrt n^*_j} \Big(
   \dfrac{\partial L_{n}(T^*_j,\widehat{\theta}_n(T^*_j))}{\partial \theta}  - \dfrac{\partial L_{n}(T^*_j,\theta^*_j)}{\partial \theta}
     \Big) \Big] .$$
 As in proof of Lemma 3 of \cite{Bardet2009}, it is now easy to show that:
  $$ \dfrac{1}{\sqrt{n^*_j}} \dfrac{\partial L_{n}(T^*_j,\theta^*_j)}{\partial \theta}  \overset{\mathcal{ D}}{\underset{n\to \infty} \longrightarrow}
   \mathcal{N}_d \big(0, G(\theta^*_j)\big)  $$
where $G(\theta^*_j)$ is given by \eqref{eq_AN}. Thus, since
${\partial
\widehat{L}_{n}(T^*_j,\widehat{\theta}_n(T^*_j))}/{\partial
\theta}=0 $, we have:
$$
\dfrac{1}{\sqrt{n^*_j}} \dfrac{\partial
L_{n}(T^*_j,\widehat{\theta}_n(T^*_j))}{\partial \theta} =
     \dfrac{1}{\sqrt{n^*_j}} \Big( \dfrac{\partial L_{n}(T^*_j,\widehat{\theta}_n(T^*_j))}{\partial
    \theta}  - \dfrac{\partial \widehat{L}_{n}(T^*_j,\widehat{\theta}_n(T^*_j))}{\partial \theta} \Big)
     \overset{\text{a.s.}}{\underset{n\to \infty} \longrightarrow} 0.
$$
We conclude using Lemma \ref{lem_conv_L} and the fact that
$1/\sqrt{n}=O(v_n/n)$. $\blacksquare$


\begin{thebibliography}{10}

\bibitem{Bai1998}
{\sc Bai J. and Perron P.}
\newblock Estimating and testing linear models with multiple structural changes.
\newblock {\em Econometrica 66\/} (1998), 47--78.




\bibitem{Bardet2009}
{\sc Bardet, J.-M. and Wintenberger, O.}
\newblock Asymptotic normality of the quasi-maximum likelihood estimator for multidimensional causal processes.
\newblock {\em Ann. Statist. 37}, (2009), 2730--2759.

\bibitem{Basseville1993}
{\sc Basseville, M. and Nikiforov, I.}
\newblock Detection of Abrupt Changes: Theory and Applications.
\newblock Prentice Hall, Englewood Cliffs, NJ, 1993.

\bibitem{Berkes2003}
{\sc Berkes, I., Horv\'ath, L., and Kokoszka, P.}
\newblock \mbox{GARCH} processes: structure and estimation.
\newblock {\em Bernoulli 9\/} (2003), 201--227.


\bibitem{Billingsley1968}
{\sc Billingsley}.
\newblock {\em Convergence of Probability Measures}.
\newblock John Wiley \& Sons Inc., New York, 1968.



\bibitem{David1995}
{\sc Davis, R. A., Huang, D. and Yao, Y.-C.}
\newblock Testing for a change in the parameter values and order of an autoregressive model.
\newblock {\em Ann. Statist. 23}, (1995), 282--304.

\bibitem{David2008}
 {\sc Davis, R. A., Lee, T. C. M. and Rodriguez-Yam, G. A.}
\newblock  Break detection for a class of nonlinear time  series models.
\newblock {\em Journal of Time Series Analysis 29},  (2008), 834--867.


\bibitem{Doukhand}
{\sc Doukhan, P., and Wintenberger, O.}
\newblock Weakly dependent chains with infinite memory.
\newblock {\em Stochastic Process. Appl. 118},  (2008) 1997-2013.


\bibitem{Duflo1990}
 {\sc Duflo, M.}
\newblock  Méthodes récursives aléatoires.
\newblock {\em Masson, Paris},  (1990); English Edition, Springer, 1996.


\bibitem{Feller1966}
{\sc Feller, W.}
\newblock {\em An Introduction to Probability Theory and its Applications},
  vol.~2.
\newblock Wiley, 1966.


\bibitem{Francq2004}
{\sc Francq, C., and Zako\"{i}an, J.-M.}
\newblock Maximum likelihood estimation of pure garch and arma-garch processes.
\newblock {\em Bernoulli 10\/} (2004), 605--637.

\bibitem{Giraitis2006}
{\sc Giraitis, L., Leipus, R., and Surgailis, D.}
\newblock Recent advances in \mbox{ARCH} modelling.
\newblock In {\em Long-Memory in Economics\/} (2006), G.~Teyssi\`ere and
  A.~Kirman, Eds., Springer Verlag., pp.~3--38.


\bibitem{Hinkley1970}
 {\sc Hinkley, D. V.}
\newblock Inference about the change-point in a sequence of random variables.
\newblock {\em Biometrika 57}, (1970), 1--17.




\bibitem{Jeantheau1998}
{\sc Jeantheau, T.}
\newblock Strong consistency of estimators for multivariate arch models.
\newblock {\em Econometric Theory 14}, 1 (1998), 70--86.

\bibitem{Leipus2000}
 {\sc Kokoszka, P. and Leipus, R.}
\newblock Change-point estimation in ARCH models.
\newblock {\em Bernoulli 6}, (2000), 513--539.

\bibitem{Kounias1969}
{\sc  Kounias E. G. and Weng T.-S.}
\newblock An inequality and almost sure convergence.
\newblock {\em Annals of Mathematical Statistics 40}, (1969), 1091--1093.


\bibitem{Lavielle2000}
{\sc Lavielle, M.  and  Ludeña, C.}
\newblock The multiple change-points problem for the spectral distribution.
\newblock {\em Bernoulli 6}, (2000) 845--869.

\bibitem{LavielleMoulines2000}
{\sc Lavielle, M. and Moulines, E.}
\newblock Least squares estimation of an unknown number of shifts in a time series.
\newblock {\em Journal of Time Series Analysis 21},  (2000) 33--59.


\bibitem{Nelson1992}
{\sc Nelson, D. B. and Cao, C. Q.}
\newblock Inequality Constraints in the Univariate GARCH Model.
\newblock {\em Journal of Business $\&$ Economic Statistics 10\/}, (1992), 229--235.


\bibitem{Page1955}
{\sc Page, E. S.}
\newblock A test for a change in a parameter occurring at an unknown point.
\newblock {\em Biometrika 42}, (1955), 523--526.

\bibitem{Page1957}
{\sc Page, E. S.}
\newblock On problems in which a change in a parameter occurs at an unknown point.
\newblock {\em Biometrika 44}, (1957), 248--252.

\bibitem{Robinson2006}
{\sc Robinson, P., and Zaffaroni, P.}
\newblock Pseudo-maximum likelihood estimation of \mbox{ARCH($\infty$)} models.
\newblock {\em Ann. Statist. 34\/} (2006), 1049--1074.

\bibitem{Straumann2006}
{\sc Straumann, D., and Mikosch, T.}
\newblock Quasi-maximum-likelihood estimation in conditionally heteroscedastic time series: A stochastic recurrence equations approach.
\newblock {\em Ann. Statist. 34}, 5 (2006), 2449--2495.


\bibitem{Yao1988}
{\sc Yao, Y.C. }
\newblock Estimating the number of change-points via Schwarz criterion.
\newblock {\em Statistics $\&$ Probability Letters 6}, (1988), 181--189.



\end{thebibliography}
\end{document}